\documentclass[11pt,english]{smfart}
 
\usepackage{ifpdf}

\usepackage{amssymb}
\usepackage{amsmath, amsthm, amsopn, amsfonts}
\ifpdf
\usepackage{graphicx}  \DeclareGraphicsExtensions{.pdf,.mps}

\else
\usepackage{amsfonts}
\usepackage{graphicx}  \DeclareGraphicsExtensions{.ps,.eps}

\fi

\usepackage{graphpap}

\def\Re{\textrm{Re}} 

\def\pasdegrille{\let\grille = \pasgrille}
\def\ecriture#1#2{\setbox1=\hbox{#1}
\dimen1= \wd1 \dimen2=\ht1 \dimen3=\dp1 \grille #2 \box1 }
\def\aat#1#2#3{
\divide \dimen1 by 48 \dimen3=\dimen1 \multiply \dimen1 by #1
\advance \dimen1 by -\dimen3 \divide \dimen1 by 101 \multiply
\dimen1 by 100 \divide \dimen2 by \count11 \multiply \dimen2 by #2
\setbox0=\hbox{#3}\ht0=0pt\dp0=0pt
   \rlap{\kern\dimen1 \vbox to0pt{\kern-\dimen2\box0\vss}}\dimen1= \wd1
\dimen2=\ht1}
\def\pasgrille{
\count12= \dimen1 \divide \count12 by 50 \divide \dimen2 by
\count12 \count11 =\dimen2 \ \divide \dimen1 by 48
\setlength{\unitlength}{\dimen1} \smash{\rlap{\ }} \dimen1= \wd1
\dimen2=\ht1 }
\def\grille{
\count12= \dimen1 \divide \count12 by 50 \divide \dimen2 by
\count12 \count11 =\dimen2 \ \divide \dimen1 by 48
\setlength{\unitlength}{\dimen1}
\smash{\rlap{\graphpaper[1](0,0)(50, \count11)}} \dimen1= \wd1
\dimen2=\ht1 }

\pasdegrille

\def\Re{\textrm{Re}} 
\def\11{{\rm 1~\hspace{-1.4ex}l} }
\def\R{\mathbb R}
\def\C{\mathbb C}
\def\Z{\mathbb Z}
\def\N{\mathbb N}

\def\T{\mathbb T}
\def\O{\mathbb O}

\newtheorem{theoreme}{Theorem}

\renewcommand{\theequation}{\arabic{equation}}

\newtheorem{proposition}{Proposition}
\newtheorem{lemme}[proposition]{Lemma}

\newtheorem{remarque}[proposition]{Remark}

\numberwithin{equation}{section}
\numberwithin{proposition}{section}
\newcommand{\rien}[1]{\relax}

\def\smfbyname{by}
\def\refname{References}
\begin{document}
\title[On the Laplace-Beltrami eigenfunctions]
{Restrictions of the Laplace-Beltrami eigenfunctions to submanifolds}
\author{N. Burq}
\address{Universit{\'e} Paris Sud,
Math{\'e}matiques,
B{\^a}t 425, 91405 Orsay Cedex et Institut Universitaire de France}
\email{Nicolas.burq@math.u-psud.fr}
\author{P. G{\'e}rard}
\address{Universit{\'e} Paris Sud,
Math{\'e}matiques,
B{\^a}t 425, 91405 Orsay Cedex}
\email{Patrick.gerard@math.u-psud.fr}
\author{N. Tzvetkov}
\address{ D\'epartement de Math\'ematiques, Universit\'e Lille I, 59 655
Villeneuve d'Ascq Cedex, France}
\email{Nikolay.tzvetkov@math.univ-lille1.fr}
\begin{abstract}
We estimates the $L^p$ norm ($2\leq p \leq +\infty$) of the restriction to a curve of the eigenfunctions of the Laplace Beltrami operator on a riemannian surface. If the curve is a geodesic, we show that on the sphere these estimates are sharp. If the curve has non vanishing geodesic curvature, we can improve our results. All our estimates are shown to be optimal for the sphere. Moreover, we sketch their extension to higher dimension
\end{abstract} 
\begin{altabstract}
On prouve une estimation de la norme $L^p$ ($2\leq p \leq + \infty$) de la restriction \`a une courbe des fonctions propres de l'op\'erateur de Laplace Beltrami sur une surface riemannienne. Si la courbe est une g\'eod\'esique de la sph\`ere, on montre que nos estimations sont optimales. En revanche, si la courbe poss\`ede une courbure g\'eod\'esique non nulle, on am\'eliore le r\'esultat. Toutes nos estim\'ees sont optimales sur la sph\`ere. Nous en esquissons par ailleurs des g\'en\'eralisations aux dimensions sup\'erieures. 
\end{altabstract}

\maketitle
\section{Introduction}
Let $(M,g)$ be a compact smooth Riemannian manifold (without boundary) of dimension
$d$. Let us denote by ${\mathbf \Delta}$ the Laplace operator associated to the metric~$g$. 
Let $\Sigma: [a,b]^k\rightarrow M$ be a smooth embedded sub-manifold of dimension $k$. The metric $g$ endows $M$ and $\Sigma$ with canonical measures  and consequently we can define the Lebesgue spaces $L^p(M)$ and $L^p(\Sigma)$, $1\leq p \leq + \infty$, of functions on $M$ and $\Sigma$ respectively. Let $(\varphi_{\lambda})$, $\lambda\geq 0$, be the eigenfunctions of ${\mathbf \Delta}$ such that $-{\mathbf \Delta}\varphi_{\lambda}=\lambda^2\varphi_{\lambda}$. This paper fits in the line of researches dealing with possible concentrations of the eigenfunctions of the Laplace operator on a manifold. There are many ways of measuring possible concentrations. One of the most popular is by describing semi-classical (Wigner) measures (see the works by Shnirelman~\cite{Sh74}, Zelditch~\cite{Ze87}, Colin de Verdi\`ere~\cite{CdV85}, G\'erard-Leichtnam~\cite{GeLe93-1}, Zelditch-Zworski~\cite{ZeZw96}, Helffer-Martinez-Robert~\cite{HeMaRo87}, Sarnak~\cite{Sa95},  Lindenstrauss~\cite{Li06} and Anantharaman~\cite{An04}). Another way is the study of the potential growth of $\|\varphi_\lambda\|_{L^p(M)}$, see the works by Sogge~\cite{So88, So93}, Sogge-Zelditch~\cite{SoZe01}, the authors~\cite{BuGeTz03-1, BuGeTz03-2, BuGeTz04}. In the present paper we propose a third way (see also the work by Reznikov~\cite{Re04}) and study the possible growth of the $L^p$ norm ($2\leq p \leq + \infty$) of the restrictions of $\varphi_\lambda$ to submanifolds of $M$. In most of this article we will concentrate to the simplest case of curves on a surface (i.e., the case when $\Sigma$ is a smooth curve
$\gamma \,: \, [a,b]\longrightarrow M$ parametrized by arc length. Our first result reads as follows. 
\begin{theoreme}\label{thm1}
There exists a constant $C$ such that for every $\varphi_{\lambda}$,
\begin{equation}\label{1}
\|\varphi_\lambda\|_{L^p( \gamma)}\leq C(1+\lambda)^{\delta(p)}
\|\varphi_{\lambda}\|_{L^2(M)}\,\, ,
\end{equation}
where 
\begin{equation} \delta(p)= 
\begin{cases} \frac 1 2 - \frac 1 p &\text{ if $4\leq p \leq + \infty$}\\
\frac 1 4 &\text{ if $2\leq p \leq 4$}.
\end{cases}
\end{equation}
Moreover (\ref{1}) is sharp if $M$ is the standard sphere $\mathbb{S}^2$
when $\gamma$ is any curve for $4\leq p \leq + \infty$ and when $\gamma$ is a geodesic curve for $2\leq p <4$.
\end{theoreme}
Notice that for example if $p=2$ then the trace theorem gives
the bound $(1+\lambda)^{1/2}$ instead of $(1+\lambda)^{1/4}$ and
thus Theorem 1 can be seen as an improvement of the trace theorem when the
traces are taken from Laplace-Beltrami eigenfunctions.

Recall that the Weyl pointwise bound provides a constant $C$ such that 
$$
\|\varphi_{\lambda}\|_{L^{\infty}(M)}\leq
C(1+\lambda)^{\frac{1}{2}}\|\varphi_{\lambda}\|_{L^2(M)}\,\, .
$$
Therefore Theorem~\ref{thm1} may also be seen as an effect of the averaging
along~$\gamma$ of $|\varphi_{\lambda}|$ which may\footnote{If for example $M$ is
the standard sphere $\mathbb{S}^2$ then $\|\varphi_{\lambda}\|_{L^{\infty}(\mathbb{S}^2)}$ may grow as $\lambda^{1/2}$
for $\lambda\rightarrow\infty$.} grow pointwisely much faster than
$\lambda^{\frac{1}{4}}$ as $\lambda\rightarrow\infty$.
\par
Let us denote by $\frac{D}{dt}$ the covariant derivative along $\gamma$ (recall that $\gamma$ is parametrized by arc length). 
For $2\leq p <4$, estimate (\ref{1}) is optimal for geodesic curves on the standard sphere $\mathbb{S}^2$,
i.e. curves such that $\frac{D}{dt}\gamma'=0$ . The second goal of this paper
is to show that estimate (\ref{1}) in the range $2\leq p <4$ can be significantly improved for curves
with non vanishing geodesic curvature. We have
the following statement.
\begin{theoreme}\label{thm2}
Let $\gamma$ be such that
\begin{equation}\label{curvature}
g\big({\frac{D}{dt}}\gamma',\frac{D}{dt}\gamma^{'}\big)\neq 0\, .
\end{equation}
There exists a constant $C$ such that for every $\varphi_{\lambda}$, $2\leq p \leq 4$,
\begin{equation}\label{1-bis}
\|\varphi_\lambda\|_{L^p(\gamma)}\leq C(1+\lambda)^{\widetilde{\delta}(p)}
\|\varphi_{\lambda}\|_{L^2(M)}\,\, ,
\end{equation}
\begin{equation} \widetilde{\delta}(p)= 
\frac 1 3 - \frac 1 {3p} 
\end{equation}
Moreover (\ref{1-bis}) is sharp in the case where $M$ is the standard sphere $\mathbb{S}^2$
and $\gamma$ is any curve with non vanishing geodesic curvature.
\end{theoreme} 
 
\begin{remarque} 
In fact both theorems above still hold if one replace the eigenfunctions by the more general spectral  $1_{ |\sqrt{- \mathbf{\Delta}} - \lambda|\leq 1}$. In that case, as will appear clearly in Section~\ref{sec.opt},  the optimality holds for such spectral projector for  {\em any compact surface}.
\end{remarque}
Let us remark that the techniques presented here extend to higher dimensions.
\begin{theoreme}\label{th.variete}
Let $(M,g)$ be a compact smooth Riemannian manifold of dimension $d$ and $\Sigma$ be a smooth submanifold of dimension $k$. There exists a constant $C>0$ such that for any $\varphi_\lambda$, we have 

\begin{equation}
\|\varphi_\lambda\|_{L^p( \Sigma)}\leq 
C(1+ \lambda)^{\rho(k,d)}\|\varphi_\lambda\|_{L^2( M)}
\end{equation}
where
\begin{equation}\label{eq.multid}
\begin{aligned} 
\rho(d-1,d)&=\begin{cases}
{\frac {d-1} 2 - \frac {d-1} p} &\text{ if } p_0= \frac {2d} {d-1} < p \leq + \infty\\
{\frac{d-1} 4 - \frac {d-2}{2p}}  &\text{ if } 2 \leq p < p_0=\frac {2d} {d-1}\\
\end{cases}\\
\rho(d-2,d)&={\frac {d-1} 2 - \frac {d-2} p}  \text{ if } 2 < p \leq + \infty\\
\rho(k,d)&={\frac {d-1} 2 - \frac {k} p}  \text{ if }  1\leq k \leq d-3
\end{aligned}
 \end{equation}
If $ p = p_0=\frac {2d} {d-1}$ and $k=d-1$, we have 
\begin{equation}
\|\varphi_\lambda\|_{L^p( \Sigma)}\leq 
C(1+ \lambda)^{\frac{d-1} {2d}}\log^{1/2}( \lambda)\|\varphi_\lambda\|_{L^2( M)}
\end{equation}
and if $ p = 2$ and $k=d-2$, we have 
\begin{equation}
\|\varphi_\lambda\|_{L^p( \Sigma)}\leq 
C(1+ \lambda)^{\frac 1 2}\log^{1/2}( \lambda)\|\varphi_\lambda\|_{L^2( M)}
\end{equation}
Moreover, all estimates are sharp, except for the log loss if
\begin{enumerate}
\item
 $k \leq d-2$,  $M$ is the standard sphere $\mathbb{S}^d$ and $\Sigma$ is any submanifold of dimension~$k$.
\item
$k={d-1}$ and $\frac {2d }{d-1}\leq p \leq + \infty$, $M$ is the standard sphere  $\mathbb{S}^d$ and $\Sigma$ is any hypersurface.
\item $k=d-1$ and $2 \leq p <\frac{2d} {d-1}$, $M$ is the standard sphere $M= \mathbb{S}^d$ and $\Sigma$ is any hypersurface containing a piece of geodesic.
\end{enumerate} 
\end{theoreme}
\begin{figure}[ht]
\label{fig:courbes}
\begin{center}\ecriture{\includegraphics[width=6cm]{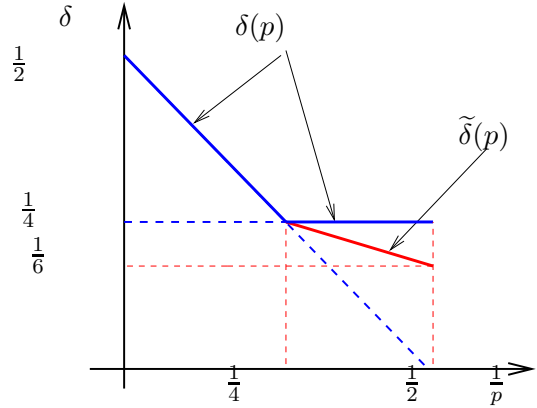}}
{
\aat{2}{35}{$\frac 1 2$}
\aat{2}{19}{$\frac 1 4$}
\aat{2}{14}{$\frac 1 6$}
\aat{22}{0}{$\frac 1 4$}
\aat{22}{39}{$\delta(p)$}
\aat{45}{27}{$\widetilde{\delta}(p)$}
\aat{38}{0}{$\frac 1 2$}
\aat{0}{40}{$\delta$}
\aat{45}{0}{$\frac 1 p$}
}
\end{center}
\caption{Curves on surfaces}
\end{figure}
\begin{figure}[ht]
\label{fig:courbebis}\begin{center}\ecriture{\includegraphics[width=6cm]{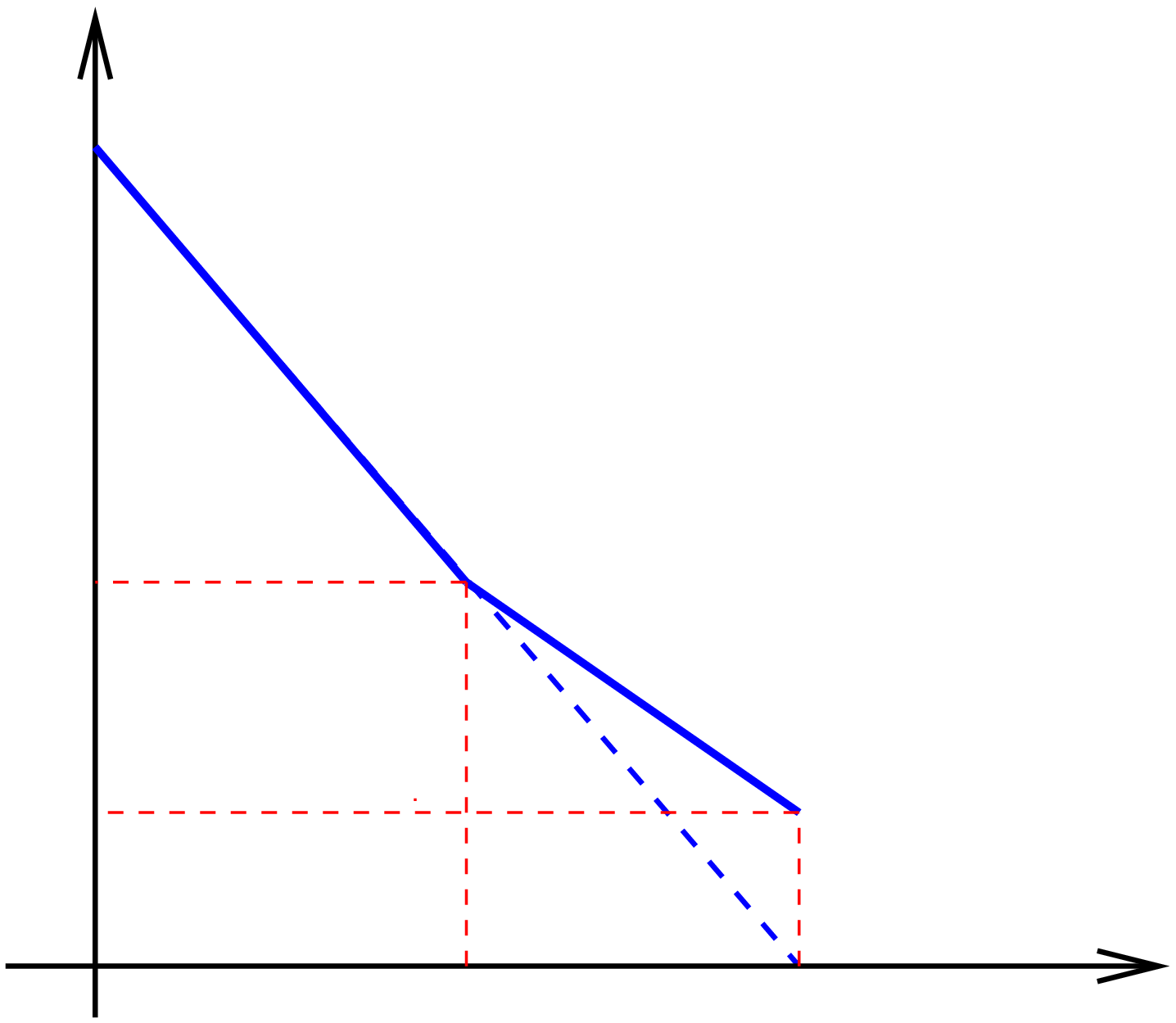} 
\includegraphics[width=6cm]{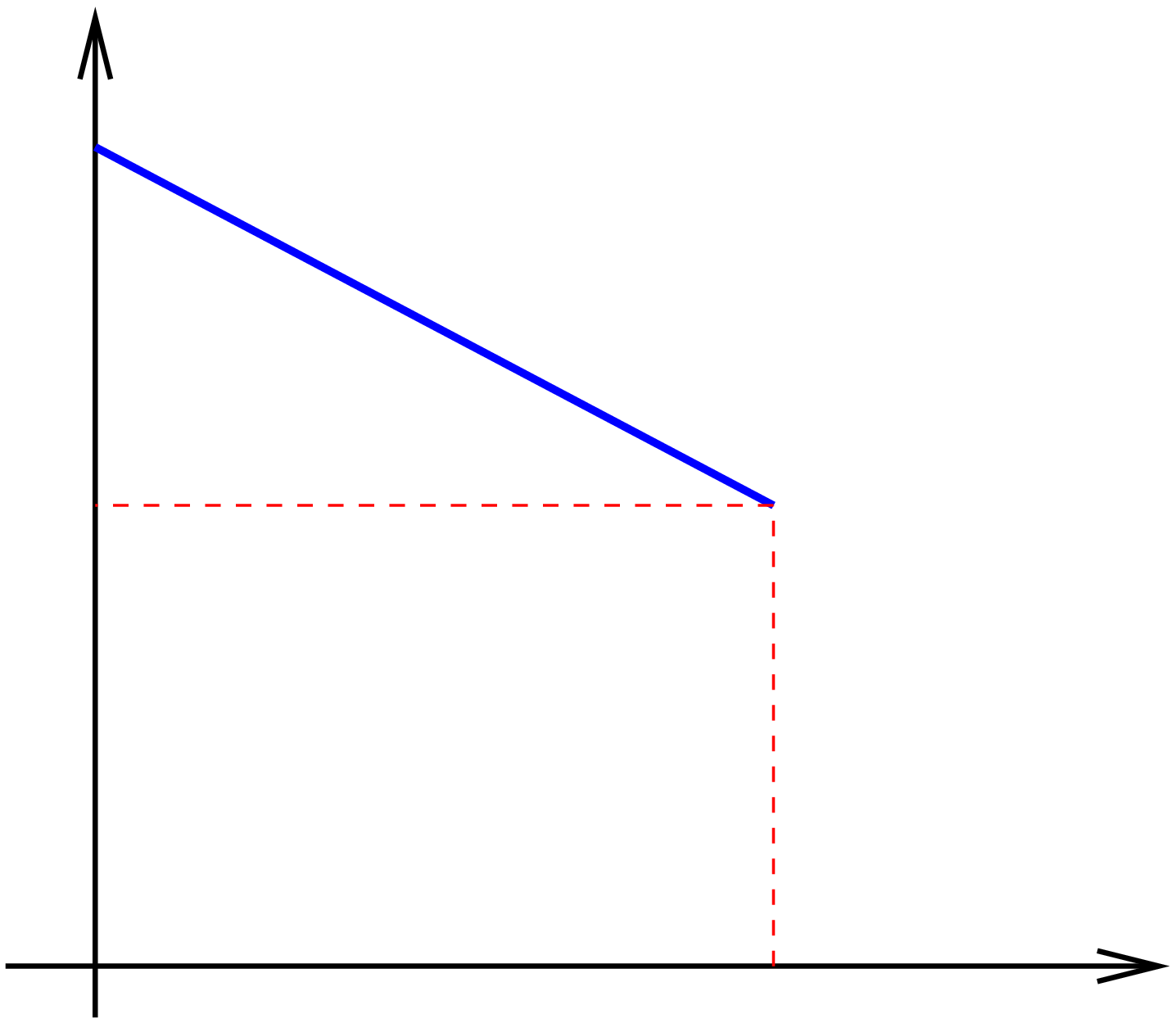}}
{\aat{3}{18}{$\frac {d-1} 2$}
\aat{11}{18}{$\mathbf{k=d-1}$}
\aat{35}{18}{$\mathbf{k\leq d-2}$}
\aat{21}{-1}{$\frac 1 2$}
\aat{50}{-1}{$\frac 1 p$}
\aat{12}{-1}{$ \frac {d-1} {2d}$}
\aat{1}{8}{$\frac{d-1} {2d}$}
\aat{1}{3}{$\frac 1 4$}
\aat{25}{18}{$\frac {d-1} 2$}
\aat{42}{-1}{$\frac 1 2$}
\aat{23}{9}{$\frac{d-1-k} 2$}
}
\end{center}
\caption{k-submanifolds in a d-manifold}
\end{figure}
\par
When  $(M,g)$ has {\em negative constant curvature} and {\em in the special case $p=2$}, estimate (\ref{1}) when $\gamma$ is a {\em geodesic} and estimate (\ref{1-bis}) for $\gamma$  a {\em geodesic circle} were recently obtained by  Reznikov~\cite{Re04}. This work was a starting point and a motivation for our present work. It was only later that we realized that in fact our results {\em in the special case of a hypersurface and $p=2$} can be deduced from Theorem 1 by Tataru~\cite{Ta97} (considering an eigenfunction as a solution of the wave equation with a trivial time dependence). Remark that in turn Tataru's result is obtained as a consequence of classical estimates on Fourier integral operators with folded canonical relations by Melrose and Taylor~\cite{Ta80-1} for the analog of our Theorem~2, (which in turn can be obtained from normal form results by Melrose~\cite{Me76}, see also H\"ormander,~\cite[Theorem 25.3.11]{Ho85-1}) and by Greenleaf and Seeger~\cite{GrSe94} for the analog of our Theorem~1. To keep our paper essentially self contained, avoid the reader to go into the normal form machinery and maintain the technical nature of the exposition at a (rather) basic level, we elected to give a complete  proof (including the case $p=2$).

It should be pointed out that estimates (\ref{1}) and (\ref{1-bis}) are far from being optimal in the case of the
flat torus $\T^2=\R^2/(2\pi \Z)^2$. In this case we have a strong improvement
of the Weyl bound on the $L^{\infty}$ norm of the eigenfunctions. More
precisely for every $\varepsilon>0$ there exists a constant $C_{\varepsilon}$
such that for every $\varphi_\lambda$,
\begin{equation}\label{2}
\|\varphi_{\lambda}\|_{L^{\infty}(\T^2)}\leq
C_{\varepsilon}(1+\lambda)^{\varepsilon}\|\varphi_{\lambda}\|_{L^2(\T^2)}\, .
\end{equation}
Indeed, in the case of the flat torus we have the explicit representation of
the eigenfunctions,
$$
\varphi_{\lambda}(x,y)=\sum_{m^2+n^2=\lambda^2}c_{m,n}e^{i(mx+ny)},\quad
c_{m,n}\in\C\, .
$$
Therefore (\ref{2}) results from the Plancherel identity 
$$\|\varphi_{\lambda}\|_{L^2(\T^2)}^{2}=\sum_{m^2+n^2=\lambda^2}|c_{m,n}|^2$$
the Cauchy-Schwarz inequality and the divisor bound 
$$
\#\big( (m,n)\in \Z^2\, : \,
m^2+n^2=\lambda^2\big)\leq C_{\varepsilon}(1+\lambda)^{\varepsilon},\quad \varepsilon>0
$$
in the ring of the Gaussian integers.
Estimate (\ref{2}) thus implies that for every curve $\gamma$ on $\T^2$, we have the bound
\begin{equation}\label{3}
\int_{\gamma}|\varphi_{\lambda}|^{2}\leq C_{\varepsilon}(1+\lambda)^{\varepsilon}
\|\varphi_{\lambda}\|_{L^2(\T^2)}^{2}\,\, .
\end{equation}
It would be interesting to decide whether one may replace
$C_{\varepsilon}(1+\lambda)^{\varepsilon}$ in (\ref{3}) by a constant uniform
with respect to $\lambda$.
\par
Let us  mention that estimate (\ref{2}) (and thus
(\ref{3})) is conjectured (see \cite{Sa95,Re04}) to hold if  the flat torus $\T^2$ is replaced by an
arbitrary negatively curved manifold. 
\par
The paper is organized as follows. In Section~\ref{sec.2} we recall the parametrix for a smoothed spectral projector (as written in~\cite{So93} for example). Then we use this representation and a classical $TT^\star$ argument to prove~\eqref{1} in the Section~\ref{sec.geod}. In the Section~\ref{sec.courbe} we precise the form of the parametrix and deduce~\eqref{1-bis} by analyzing precisely the oscillations which appear in the phase. In the Section~\ref{sec.opt}, this precise description is used to prove the optimality of~\eqref{1} and ~\eqref{1-bis} on the sphere.  In the last section we show how the two dimensional techniques can be adapted to prove Theorem~\ref{th.variete}.\par
{\bf Acknowledgements.} A previous version of this article was only dealing with curves on surfaces. We thank one referee for pointing to us the question of higher dimension submanifolds and another one for various comments which improved the presentation of the paper and for proposing a counter example in section~\ref{sec.opt}. We also thank Melissa Tacy for pointing out an incorrect point in a previous proof of Lemma~\ref{lem.optbis}

\section{Representation of $\varphi_{\lambda}$ in local coordinates}\label{sec.2}
We have the following representation of $\varphi_{\lambda}$ in local coordinates.
\begin{theoreme}\label{sogge} There exits a function $\chi \in \mathcal{S}(\mathbb{R})$ such that $\chi(0)=1$ and for any $x_0\in M$ there exists systems of coordinates near $x_0$,
$$
W\subset V=\{x\in\R^d\, : \, x\sim 0\},
$$
a smooth function
$$ a\,:\, W_x\times V_y\times \mathbb{R}_\lambda^+\longrightarrow\C $$
supported in the set 
\begin{equation*}
\{(x,y)\in W\times V\, : \, |x| \leq c_0\varepsilon \leq c_1 \varepsilon\leq |y|\leq
c_2\varepsilon\ll 1 \}\, 
\end{equation*}and such that
$$ \forall \alpha \in \mathbb{N}^{2d}, \exists C>0; \forall \lambda \geq 0,  |\partial_{x,y}a(x,y, \lambda)|\leq C$$
and an operator 
$$
\mathcal{R}_\lambda\,:\, L^2(W)\longrightarrow L^2(W)
$$
such that, $$
\|\mathcal{R}_\lambda f\|_{L^\infty(W)}\leq C \|f\|_{L^2(W)}
$$ and with $U:=W\cap \{ x\,:\, |x|\leq c\,\varepsilon\}$, for every $x\in U$,
\begin{equation}\label{eq.descrip}
\chi_\lambda(f):=\chi( \sqrt{- \mathbf{\Delta} }- \lambda) f = \lambda^{\frac{d-1}{2}}\int_{y\in
  V}e^{i\lambda\psi(x,y)}a(x,y, \lambda)
f(y)dy+\mathcal{R}_\lambda (f)\, ,
\end{equation}
where
$\psi(x,y)=-d_{g}(x,y)$ is the geodesic distance with respect to $g$ between
$x$ and $y$. Furthermore, the symbol $a(x,y, \lambda)$ is real non negative and does not vanish for $|x|\leq c\varepsilon$ and $ d_{g}(x,y) \in [c_3 \varepsilon,  c_4 \varepsilon]$. 
\end{theoreme}
Remark that since $\chi( \sqrt{- \mathbf{\Delta} }- \lambda) \varphi_\lambda = \varphi_\lambda$, relations~\eqref{1} and~\eqref{1-bis} are consequences of similar estimates on the norm of the operator $\chi( \sqrt{- \mathbf{\Delta} }- \lambda) $ from $L^2(M)$ to $L^p( \gamma)$. Remark also that, since $\mathcal{R}_\lambda$ satisfies better estimates than we claimed in Theorems~1 and~2, we will restrict the study to the analysis of the principal part in~\eqref{eq.descrip}, $\mathcal{T}_\lambda$
$$
\mathcal{T}_{\lambda}(f)(x):=\int_{y\in  V}e^{i\lambda\psi(x,y)}a(x,y)
f(y)dy\, .
$$\par
The proof of Theorem~\ref{sogge} may be found in \cite[Chap. 5]{So93} (except for the non vanishing of the symbol which comes from the explicit calculations made there) and is based on a
representation of $\varphi_{\lambda}$ by the aid of the propagator of the wave
equation
$$
(\partial_{t}^{2}-{\mathbf \Delta})u=0\, .
$$
More precisely if $\chi\in {\mathcal S}(\R)$, $\chi(0)=1$ is a Schwartz function such
that ${\rm supp}(\widehat{\chi})\subset[\varepsilon, 2\varepsilon]$
where $\varepsilon\ll 1$, then we can write 
\begin{equation*}
\chi(\sqrt{-{\mathbf \Delta}}-\lambda)f
 = 
\frac{1}{2\pi}\int_{\varepsilon}^{2\varepsilon}
\Big(
e^{i\tau\sqrt{-{\mathbf \Delta}}}f
\Big)e^{-i\tau\lambda}\widehat{\chi}(\tau)d\tau\,.
\end{equation*}
For $\varepsilon\ll 1$ and $|\tau|\leq 2\varepsilon$, 
introducing local coordinates near $x_0\in M$, 
we can represent $e^{i\tau\sqrt{-{\mathbf \Delta}}}$ as a Fourier integral operator
(see e.g. \cite{Ho68}). A stationary phase
argument (see \cite[Chap. 5]{So93}) in this representation then achieves the
proof of Theorem~\ref{sogge}.
\section{Proof of (\ref{1})}\label{sec.geod}
In this section we give the proof of (\ref{1}). 
Let $x_0\in \gamma$. Using a partition of unity, we may assume that $\gamma$ is contained in the domain $U$ of our coordinate patch defined in Theorem~\ref{sogge} and $x_0=0$. We assume that $\gamma$ is parametrized by arc length, $s$ and define 
$$ T_\lambda(f) (s) = \mathcal{T}_\lambda(f) (x(s)).$$
Therefore it suffices to prove that
\begin{equation}\label{2.1}
\|T_{\lambda}f\|_{L^p(\gamma)}
\leq
 C\lambda^{\delta(p)-\frac{1}{2}}
\left(\int_{y\in V}|f(y)|^{2}dy\right)^{1/2}\, .
\end{equation}
We represent $y$ in geodesic polar coordinates as 
$y=\exp_{0}(r\omega)$, $r >0$, $\omega \in \mathbb{S}^{1}$ and $c_1\varepsilon\leq
r\leq c_2\varepsilon$. We can write
$$
(T_{\lambda}f)(x)=\int_{c_1\varepsilon}^{c_2 \varepsilon}
(T_{\lambda}^{r}f_r)(x)dr
$$
where
$$
(T_{\lambda}^{r}f)(x)=\int_{\mathbb{S}^1}
e^{i\lambda\psi_{r}(x,\omega)}a_r(x,\omega)
f(\omega)d\omega
$$
with
$$
\psi_{r}(x,\omega)=\psi(x,y),\quad f_{r}(\omega)=f(y),\quad 
a_r(x,\omega)=\kappa(r, \omega)a(x,y)
$$
for some smooth function $\kappa$. It suffices consequently to show the bound
\begin{equation}\label{2.2}
\|T_{\lambda}^{r}f\|_{L^p(\gamma)}
\leq C\lambda^{\delta(p)- \frac 1 2 }\left(\int_{\mathbb{S}^1}|f(\omega)|^{2}d\omega\right)^{1/2}\, .
\end{equation}
Indeed if we suppose that (\ref{2.2}) is true, then we can write
\begin{multline}
\|T_{\lambda}f\|_{L^p(\gamma)}\leq
\int_{c_1\varepsilon}^{c_2 \varepsilon}
\|T_{\lambda}^{r}f_r\|_{L^p(\gamma)}dr
\\
\leq
C\lambda^{\delta(p)-\frac{1}{2}}\int_{c_1\varepsilon}^{c_2
  \varepsilon}\|f_r\|_{L^2(\mathbb{S}^1)}dr\leq c\lambda^{\delta(p)-\frac{1}{2}}
\|f\|_{L^2(V)}
\end{multline}
which is (\ref{2.1}). Thus it remains to prove (\ref{2.2}). For that purpose
we use a duality argument. Computing $\|(T_{\lambda}^{r})^{\star}f\|_{L^2}^2$, we can apply the usual $T$-$T^\star$ argument:
\begin{equation}\label{eq.3.4}
\begin{aligned}
\|T_{\lambda}^{r}\|^2_{\mathcal{L} (L^2(\mathbb{S}^1);L^{p}(\gamma))}&=\|(T_{\lambda}^{r})^{\star}\|^2_{\mathcal{L} (L^{p'}(\gamma);L^2(\mathbb{S}^1))}\\
&=\|T_{\lambda}^{r}(T_{\lambda}^{r})^{\star}\|_{\mathcal{L} (L^{p'}(\gamma); L^{p}(\gamma))}
\end{aligned}
\end{equation}
A direct computation shows that
$$
(T_{\lambda}^{r}(T_{\lambda}^{r})^{\star}f)(x(t))=
\int_{a}^{b}K(t,\tau)f(x(\tau))d\tau\, ,
$$
where
$$
K(t,\tau)=
\int_{\mathbb{S}^1}e^{i\lambda[\psi_{r}(x(t),\omega)-\psi_{r}(x(\tau),\omega)]}
a_{r}(x(t),\omega)\overline{a_r}(x(\tau),\omega)
d\omega\, .
$$
One can calculate the traces at zero of the first order derivatives of the map
$$
x\longrightarrow \psi_{r}(x,\omega)\, .
$$
More precisely, we have the following statement.
\begin{lemme}\label{jet}
Suppose that $g(0)={\rm Id}$. Then for $\varepsilon\ll 1$, 
$$
(\partial_{x_1}\psi_r(0,\omega),\partial_{x_2}\psi_r(0,\omega))=\omega,\quad
\omega=(\omega_1,\omega_2)\in \mathbb{S}^1\, .
$$
\end{lemme}
\begin{proof}
The proof of Lemma~\ref{jet} may be found in \cite{BuGeTz03-1}. For the sake of completeness, we recall it below.
Let $y=\exp_{0}(r \omega)$ and  $u=u(x,y)\in T_{y}M$ be the unit vector such that
$$
\exp_{y}(-\psi(x,y)u(x,y))=x.
$$
Differentiating with respect to $x$ this identity, we get for $x=0$, and any $h\in T_{0}M$,
\begin{equation}\label{eq4.100}
h=
T_{ru(0,y)}(\exp_{y})
\left[
-d_{x}\psi(0,y)\cdot hu(0,y)-r\, T_{x}u(0,y)\cdot h
\right].
\end{equation}
But since $y= \exp_{0}(r \omega)$
\begin{equation}
T_{ru(0,y)}(\exp_{y})\cdot u(0,y)=-\omega
\end{equation}
Take the scalar product with $\omega$ in~\eqref{eq4.100}. Using Gauss' Lemma (see for example~\cite[3.70]{GaHuLa90}), we get
$$
d_{x}\psi(0,y)\cdot h=g_{0}(h,\omega),
$$
i.e.
$$
\nabla_{x}\psi(0,y)=\omega.
$$
\end{proof}
Note that by performing a linear change of coordinates, we can assume that $g(0)={\rm Id}$.
Using Lemma~\ref{jet}, we can estimate the kernel $K(t,\tau)$.
\begin{lemme}\label{phase-st}
There exists $\delta>0$ such that for $|t-\tau|<\delta$, 
$$|K(t,\tau)|\lesssim (1+\lambda|t-\tau|)^{-\frac{1}{2}}.$$
\end{lemme}
\begin{proof}
For $(x,x')\in U\times U$, we set
$$
K(x,x')=\int_{\mathbb{S}^1}e^{i\lambda[\psi_{r}(x,\omega)-\psi_{r}(x',\omega)]}
a_{r}(x,\omega)\overline{a_r}(x',\omega)
d\omega\, .
$$
We are going to show that for $x, x'$ close to zero, 
\begin{equation}\label{disp}
|K(x,x')|\lesssim (1+\lambda|x-x'|)^{-\frac{1}{2}}
\end{equation}
which, in view of the property $g(\dot{x}(t),\dot{x}(t))=1$, implies the assertion of the lemma.
Let us now give the proof of (\ref{disp}). We parametrize the circle $\mathbb{S}^1$ as
$$
\omega=\omega(w)=(\cos w,\sin w),\quad w\in[0,2\pi]\, .
$$
The Taylor formula allows to write
\begin{equation}\label{edno}
\psi_{r}(x,\omega)-\psi_{r}(x',\omega)=\langle x-x',\Psi(x,x',\omega)\rangle,
\end{equation}
where
\begin{equation}\label{dve}
\Psi(x,x',\omega)
=\int_{0}^{1}\nabla_{x}\psi_{r}(x'+\theta(x-x'),\omega)d\theta\, .
\end{equation}
In (\ref{edno}), (\ref{dve}), $\langle\cdot,\cdot\rangle$ stands for the $\R^2$
scalar product and $\nabla_{x}=(\partial_{x_1},\partial_{x_2})$.
We can therefore write
$$
\psi_{r}(x,\omega)-\psi_{r}(x',\omega)=|x-x'|\Phi(x,x',\sigma,\omega),
$$
where
$$
\Phi(x,x',\sigma,\omega)
=
\langle\sigma,\Psi(x,x',\omega)\rangle,\quad
\sigma=\frac{x-x'}{|x-x'|}\in \mathbb{S}^1\, .
$$
For $(x,x',\sigma)\in U\times U\times \mathbb{S}^1$ we set
$$
K_{1}(x,x',\sigma):=\int_{\mathbb{S}^1}e^{i\lambda\Phi(x,x',\sigma,\omega)}
a_{r}(x,\omega)\overline{a_r}(x',\omega)
d\omega\, .
$$
It suffices therefore to show that there exist $\varepsilon\ll 1$ and $C>0$
such that for every $(x,x',\sigma)\in U\times U\times \mathbb{S}^1$, every $\lambda\geq 0$,
$$
|K_1(x,x',\sigma)|\leq C(1+\lambda)^{-\frac{1}{2}}\, .
$$
Suppose that
$$
\sigma=(\cos\alpha,\sin\alpha),\quad \alpha\in[0,2\pi]\, .
$$
Thanks to Lemma~\ref{jet}, 
$$
\Phi(0,0,\sigma,\omega(w))=\langle\sigma,\omega\rangle=\cos(w-\alpha)\,.
$$
Next, we can write
\begin{equation}\label{parvo}
\Phi'_{w}(0,0,\sigma,\omega(w))=-\sin(w-\alpha)
\end{equation}
and
\begin{equation}\label{vtoro}
\Phi''_{ww}(0,0,\sigma,\omega(w))=-\cos(w-\alpha) \, .
\end{equation}
The main point is that $\Phi'_{w}(0,0,\sigma,\omega(w))$ and
$\Phi''_{ww}(0,0,\sigma,\omega(w))$ cannot vanish simultaneously. 
In view of (\ref{parvo}) and (\ref{vtoro}) 
we can represent the circle $\mathbb{S}^1$ as a disjoint union of $4$ segments
$I_1(\sigma)$, $I_2(\sigma)$, $I_3(\sigma)$ and $I_4(\sigma)$, where $I_2(\sigma)$ and $I_4(\sigma)$ are neighborhoods
of $\alpha$ and $\alpha+\pi$ respectively, so that
$$
|\Phi'_{w}(0,0,\sigma,\omega(w))|\geq \frac{\sqrt{2}}{2}\, ,\quad w\in I_1(\sigma)\cup I_3(\sigma) 
$$
and
$$
|\Phi''_{ww}(0,0,\sigma,\omega(w))|\geq \frac{\sqrt{2}}{2}\,,\quad w\in I_2(\sigma)\cup I_4(\sigma).
$$
Therefore, by continuity, there exists $\varepsilon\ll 1$ such that for every
$\sigma\in \mathbb{S}^1$, if 
$$
(x,x')\in U\times U
$$ 
then
\begin{equation}\label{VdC1}
|\Phi'_{w}(x,x',\sigma,\omega(w))|\geq \frac{\sqrt{2}}{4}  \, ,\quad w\in I_1(\sigma)\cup I_3(\sigma)
\end{equation}
and
\begin{equation}\label{VdC2}
|\Phi''_{ww}(x,x',\sigma,\omega(w))|\geq \frac{\sqrt{2}}{4}  \,,\quad w\in I_2(\sigma)\cup I_4(\sigma).
\end{equation}
We now split the integral defining $K_1(x,x',\sigma)$ in $4$ parts according
to the splitting of $\mathbb{S}^1$ to $I_1(\sigma)$, $I_2(\sigma)$, $I_3(\sigma)$ and $I_4(\sigma)$.
Thanks to (\ref{VdC1}) and an integration by parts the contributions of $I_1(\sigma)$ and $I_3(\sigma)$ are bounded by
$$
C(1+\lambda)^{-1}\, .
$$
Thanks to (\ref{VdC2}) and the basic van der Corput lemma (see \cite[Chap VIII.1]{St93}) the
contributions of $I_2(\sigma)$ and $I_4(\sigma)$ are bounded by
$$
C(1+\lambda)^{-\frac{1}{2}}\, .
$$
This completes the proof of  Lemma~\ref{phase-st}.
\end{proof}
Let us now complete the proof of Theorem~\ref{thm1}.  We first recall the following classical consequence of Young's and Hardy-Littlewood's inequalities
\begin{lemme}\label{lem.convol}
 Consider $2\leq p \leq+\infty$ and $T$ a convolution operator
$$ Tf(x)= \int_{-\infty}^{+\infty} K(x-y) f(y) dy$$
with $K\in L^1_{\text{loc}}$. Then the norm of $T$ as an operator from $L^{p'}( \mathbb{R})$ to $L^p(\mathbb{R})$ is bounded by $\|K\|_{L^{\frac p 2}}$.
Furthermore, if $K(x)= \frac 1 {|x|^{2/p}}$ and $p \neq 2$, then the operator $T$ is still bounded from $L^{p'}( \mathbb{R})$ to $L^p( \mathbb{R})$ (despite the logarithmic divergence of $\|K\|_{L^{\frac p 2 }( \mathbb{R})}$).
\end{lemme}
Using  Lemma~\ref{phase-st}, we can write
\begin{eqnarray*}
\|T_{\lambda}^{r}(T_{\lambda}^{r})^{\star}f\|_{L^p(\gamma)}
& \lesssim &
\|\int_{a}^{b}(1+\lambda|t-\tau|)^{-\frac{1}{2}}|f(x(\tau))|d\tau\|_{L^{p}(a,b)}
\\
& \lesssim &
\|(1+\lambda|t|)^{-\frac{1}{2}}\|_{L^{\frac p 2}(a,b)}\|f\|_{L^{p'}(\gamma)}\,\, .
\end{eqnarray*}
But
$$
\|(1+\lambda|t|)^{-\frac 1  2}\|_{L^{\frac p 2}(a,b)}\leq
\lambda^{-\frac 2 p}(\int_{0}^{C\lambda}(1+\tau)^{-\frac{p}{4}}d\tau)^{\frac 2 p}
\lesssim \begin{cases}
\lambda^{-\frac{1}{2}} & \text{if $2\leq p < 4$}\\
\lambda^{-\frac{2}{p}} & \text{if $4 <p$}
\end{cases}
$$ 
Therefore for $p\neq 4$ $T_{\lambda}^{r}(T_{\lambda}^{r})^{\star}$ sends
$L^{p'}(\gamma)$ to $L^{p}(\gamma)$ with operator norm
$\leq C\lambda^{\delta(p) - \frac 1 2}$. The case $p=4$ follows from the additional property in Lemma~\ref{lem.convol}. This ends the proof of Theorem~\ref{thm1} (except for the optimality which will be dealt with in section~\ref{sec.opt}).
\qed
\section{Proof of (\ref{1-bis})}\label{sec.courbe}
This section is devoted to the proof of (\ref{1-bis}). For that purpose we need a more precise description of the kernel of the operators which appeared in the previous section. Notice that this section contains all the results in section~\ref{sec.geod}, but we preferred to keep the technical level of the exposition as low as possible for the proof of Theorem~1. To prove (\ref{1-bis}) we only need to consider the case $p=2$ (all the other cases being deduced from that one and the special case $p=4$ in Theorem~1 by interpolation).

As in the proof of
(\ref{1}), for $x_0\in\gamma$ and $\lambda\geq 1$, we consider the map
$$
(T_{\lambda}f)(x):=\int_{y\in  V}e^{i\lambda\psi(x,y)}a(x,y)
f(y)dy,\quad x\in U,
$$
where $U$ and $V$ are coordinate systems near $x_0$ given in
Theorem~\ref{sogge}, $a(x,y)$ is the amplitude introduced in the
statement of Theorem~\ref{sogge} and $-\psi(x,y)$ is the geodesic distance
between $x$ and $y$. Estimate (\ref{1-bis}) is therefore a consequence of the
estimate
\begin{equation*}
\|T_{\lambda}f\|_{L^2(\gamma)}^{2}
\leq
 C\lambda^{-\frac{2}{3}}
\int_{y\in V}|f(y)|^{2}dy\, .
\end{equation*}
Recall that the $L^{2}(\gamma)$ norm measures the contribution of $U$ to the
$L^2$ norm along $\gamma$. By duality it suffices to show that $T_{\lambda}^{\star}$ sends $L^{2}(V)$ to
$L^{2}(\gamma)$ with norm $\leq C\lambda^{-\frac{1}{3}}$. Finally, both
previous statements are equivalent to the continuity of $T_{\lambda}T_{\lambda}^{\star}$ from $L^{2}(\gamma)$ to
$L^{2}(\gamma)$ with operator norm $\leq C\lambda^{-\frac{2}{3}}$.
As in the proof of (\ref{1}) the main point in the proof of (\ref{1-bis}) is a
precise description of the kernel of $T_{\lambda}T_{\lambda}^{\star}$.
The new fact is that
an oscillatory factor in the kernel of $T_{\lambda}T_{\lambda}^{\star}$ will
be crucial to achieve the needed bound. Recall that in the proof of (\ref{1})
only pointwise bounds on the kernel of $T_{\lambda}T_{\lambda}^{\star}$ were used.
We can write 
$$
(T_{\lambda} T_{\lambda}^{\star})(f)(x(t))=
\int_{a}^{b}K(t,\tau)f(x(\tau))d\tau\, ,
$$
where
$$
K(t,\tau)={\mathcal K}(x(t),x(\tau))
$$
with
$$
{\mathcal K}(x,x')
=
\int_{V}e^{i\lambda[\psi(x,y)-\psi(x',y)]}
a(x,y)\overline{a}(x',y)  
dy\, .
$$ 
We now give a precise description of this kernel.
\begin{lemme}\label{geom}
There exist $\varepsilon\ll 1$, 
$a^{\pm},b\in C^{\infty}(\R^2\times\R^2\times \R)$
such that for $|x-x'|\gtrsim \lambda^{-1}$, $x,x'\in U$,
$$
{\mathcal K}(x,x')=\sum_{\pm}\frac{e^{\pm i\lambda \psi(x,x')}}
{(\lambda|x-x'|)^{\frac{1}{2}}}\,
a^{\pm}(x',x-x')+b(x,x',\lambda^{-1}|x-x'|^{-1})\, .
$$
Moreover $a^{\pm}$ are real, have supports of size ${\mathcal O}(\varepsilon)$ with
respect to the $x$ and $x'$ variables and $a^\pm(0,0)\geq c>0$. Finally
$$
|b(x,x',\lambda^{-1}|x-x'|^{-1})|\leq C(\lambda|x-x'|)^{-\frac 3 2}\,\, .
$$
\end{lemme}
\begin{remarque} We have $\chi_\lambda \chi_\lambda^\star= (|\chi|^2)_\lambda$. As a consequence, the lemma above amounts to describing the kernel of $|\chi|^2_\lambda$. The main difference with respect to Sogge's result in Theorem~\ref{sogge} is that the Fourier transform of $|\chi|^2$ no longer vanishes near $t=0$ which makes the description more difficult.
\end{remarque}
\begin{proof}
We represent $y$ in geodesic polar coordinates with origin at $x'$ as
$$
y=\exp_{x'}(r\omega),\quad r>0,\quad \omega\in \mathbb{S}^1\,\, .
$$
Denote by 
$$b_r(x,x', \omega)= a(x, \exp_{x'}(r\omega)) \overline{a} (x', \exp_{x'}(r\omega)) \kappa (x', r, \omega).$$
We can therefore write
\begin{equation}\label{tos}
{\mathcal K}(x,x')=
\int_{r= c_1\varepsilon}^{c_2\varepsilon}\int_{\mathbb{S}^1}
e^{i\lambda[\psi_{r}(x,\omega)-\psi_{r}(x',\omega)]}
b_r(x,x',\omega)d\omega dr \, ,
\end{equation}
where $\kappa(r)$ is a smooth positive function and 
$$
\psi_{r}(x,\omega)=\psi(x,y)=-d_{g}\big(x,\exp_{x'}(r\omega)\big)\, .
$$
We are going to evaluate the integral over $\omega$ in (\ref{tos}) by the
stationary phase formula. 
Assume that
\begin{equation*}
\omega=(\cos w,\sin w),\quad w\in[0,2\pi]\,\, .
\end{equation*}

We are going to show that for $x\neq x'$, close to each other, the map
\begin{equation}\label{map}
\omega\longrightarrow 
\psi_{r}(x,\omega)-\psi_{r}(x',\omega)
\end{equation}
from $\mathbb{S}^1$ to $\R$ has exactly two non degenerate critical points.
By the fundamental theorem of calculus
$$
\psi_{r}(x,\omega)-\psi_{r}(x',\omega)
=|x-x'|\int_{0}^{1}\langle
\sigma,\nabla_{x}\psi_{r}(x'+\theta(x-x'),\omega)\rangle d\theta\,,\quad
\sigma=\frac{x-x'}{|x-x'|}\,.
$$
Further, we consider the real valued function $\Phi_{r}$ 
on $\mathbb{S}^1\times U\times U\times \mathbb{S}^1$, defined by
$$\Phi_{r}(\sigma,x,x',\omega):=
\int_{0}^{1}\langle
\sigma,\nabla_{x}\psi_{r}(x'+\theta(x-x'),\omega)\rangle d\theta\,.
$$
Let us show that for fixed $\sigma\in \mathbb{S}^1$ and $x,x'$ close to each other, the map
\begin{equation}\label{function}
w\longrightarrow \Phi_{r}(\sigma,x,x',\omega(w))
\end{equation}
from $[0,2\pi]$ to $\R$ has exactly two non degenerate critical points.
We can assume that $g(x')={\rm Id}$ and using Lemma~\ref{jet}, we obtain 
$$
\Phi_{r}(\sigma,x',x',\omega)=\langle\sigma,\omega\rangle
$$
and therefore
\begin{equation}\label{hess}
\frac{\partial \Phi_{r}}{\partial w}(\sigma,x',x',\omega(w))
=\langle\sigma,\dot{\omega}\rangle\, ,
\end{equation}
where $\dot{\omega}=(-\sin w,\cos w)$ denotes the first derivative of $\omega$ with respect to $w$.
Hence, for fixed $\sigma\in \mathbb{S}^1$ the equation 
$$
\frac{\partial \Phi_{r}}{\partial w}(\sigma,x',x',\omega)=0
$$
has exactly two solutions $\omega=\pm \sigma$. Clearly
$$
\frac{\partial^{2} \Phi_{r}}{\partial w^2}
(\sigma,x',x',\omega)=\langle\sigma,\ddot{\omega}\rangle\,,
$$
where $\ddot{\omega}=(-\cos w,-\sin w)$ denotes the second derivative of $\omega$ with respect to
$w$. Therefore
\begin{equation*}
\Big|\frac{\partial^{2} \Phi_{r}}{\partial w^2}
(\sigma,x',x',\pm \sigma)\Big|=1 \,\, .
\end{equation*}
Thanks to the implicit function theorem\footnote{Notice the uniform dependence
of $\delta,\delta'$ with respect to $\sigma\in \mathbb{S}^1$.}, we obtain that there
exist $\delta>0$, $\delta'>0$ such that for every $\sigma\in \mathbb{S}^{1}$, 
every $|x-x'|<\delta$ there exists a unique $\omega\in \mathbb{S}^1$ satisfying
$|\omega-\sigma|<\delta'$ which is a non degenerate critical point of (\ref{function}).
Moreover there exists a unique $\omega\in \mathbb{S}^1$ satisfying
$|\omega+\sigma|<\delta'$ which is a non degenerate critical point of
(\ref{function}). Moreover for $\omega$ outside the intervals
$|\omega-\sigma|<\delta'$ and $|\omega+\sigma|<\delta'$, thanks to (\ref{hess}),
$$
\Big|\frac{\partial \Phi_{r}}{\partial w}(\sigma,x,x',\omega)\Big|\geq c>0\,,
$$
provided $x,x'$ being close to each other.
We therefore obtain that (\ref{function}) has exactly two non degenerate
critical points.
This in turn implies that for $x\neq x'$ the map (\ref{map}) has also 
exactly two non degenerate
critical points.
\par
Our next step is to construct explicitly the critical points of (\ref{map}).
Let us denote by $\Sigma$ the circle centered at $x'$ and of radius $r$ in the
geodesic coordinate system with center $x'$ (see Figure~2). The geodesic
$l$ joining 
the (different) points $x$ and $x'$ is a straight line in this coordinate system. Let us denote by $\exp_{x'}(r\omega^{\star})$ and
$\exp_{x'}(r(\omega^{\star}+\pi))$ the intersections of $l$ and $\Sigma$. 
We can clearly suppose that $x$ lies on the segment between $x'$ and
$\exp_{x'}(r\omega^{\star})$.
\begin{figure}[ht]
\label{fig.geod}
$$\ecriture{\includegraphics[width=5cm]{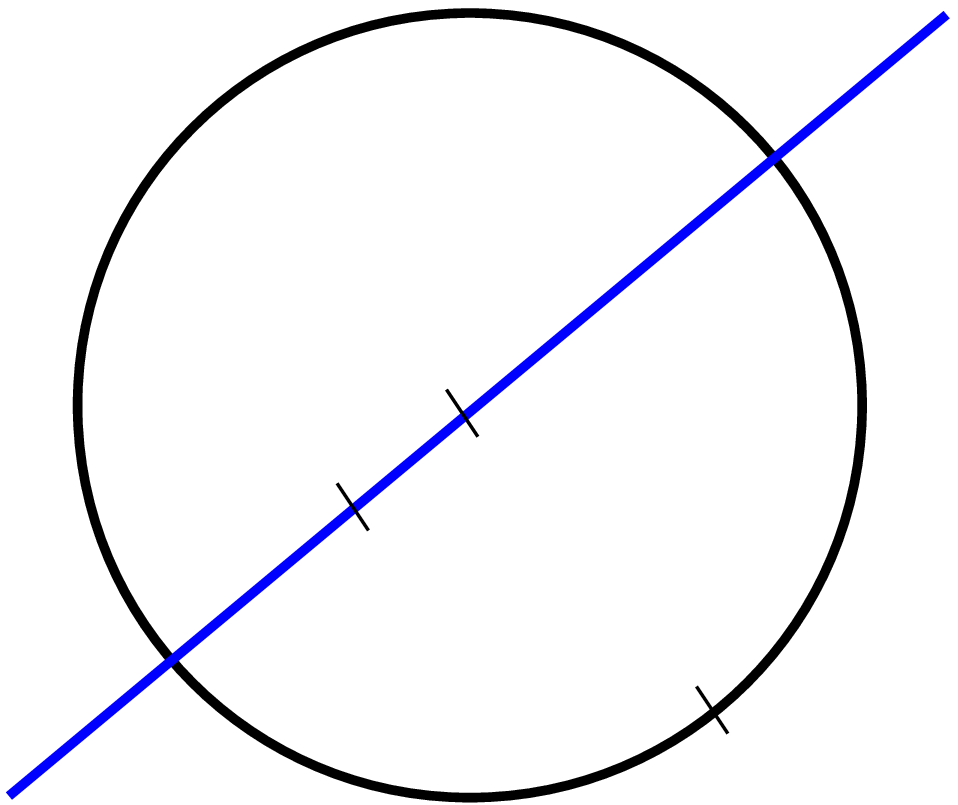}}
{\aat{45}{41}{$l$}\aat{22}{22}{$x'$}\aat{17}{18}{$x$}\aat{3}{35}{$\Sigma$}\aat{43}{32}{$\exp_{x'}(r(\omega^*+ \pi))$}\aat{42}{5}{$\exp_{x'}(r\omega)$}\aat{2}{-2}{$\exp_{x'}(r\omega^*)$}}$$ \caption{}
\end{figure} 
We are now going to show that $\omega^{\star}$
and $\omega^{\star}+\pi$ are (the) critical points of $(\ref{map})$.
Let $\exp_{x'}(r\omega)$, $\omega\in \mathbb{S}^1$ be an arbitrary point on $\Sigma$.
The triangle inequality yields
\begin{equation}\label{triangle}
d_{g}(x',\exp_{x'}(r\omega))-d_{g}(x,x')
\leq
d_{g}(x,\exp_{x'}(r\omega))\leq d_{g}(x,x')+d_{g}(x',\exp_{x'}(r\omega))\, .
\end{equation}
Notice that
$
d_{g}(x',\exp_{x'}(r\omega))=r,
$
and moreover using that  the line joining $x$ and $x'$ is a geodesic, we infer
the identities
$$
d_{g}(x,x')+d_{g}(x,\exp_{x'}(r\omega^{\star}))=r,\quad
d_{g}(x,\exp_{x'}(r(\omega^{\star}+\pi)))-d_{g}(x,x')=r\,\, .
$$
Therefore, in view of (\ref{triangle}),
\begin{equation}\label{rav}
d_{g}(x,\exp_{x'}(r\omega^{\star}))\leq d_{g}(x,\exp_{x'}(r\omega))\leq
d_{g}(x,\exp_{x'}(r(\omega^{\star}+\pi)))
\end{equation}
and equalities are possible if and only if $\omega=\omega^{\star}$ or
$\omega=\omega^{\star}+\pi$.
But (\ref{rav}) can be rewritten as
$$
\psi_{r}(x,\omega^{\star})\leq \psi_{r}(x,\omega)\leq
\psi_{r}(x,\omega^{\star}+\pi),\quad \forall\, \omega\in \mathbb{S}^1\,.
$$
On the other hand 
$$
\psi_{r}(x',\omega)=r,\quad \forall\, \omega\in \mathbb{S}^1
$$
and therefore $\omega^{\star}$ is the unique global minimum of (\ref{map})
while $\omega^{\star}+\pi$ is the unique global maximum of (\ref{map}). 
Hence, for $x\neq x'$, the critical points of (\ref{map}) are $\omega^{\star}$ and $\omega^{\star}+\pi$.
\par
We next evaluate the value of $\psi_{r}(x,\omega)-\psi_{r}(x',\omega)$ at the
critical points. Since the line joining $x$ and $x'$ is a geodesic, we obtain
that
$$
\psi_{r}(x,\omega^{\star})-\psi_{r}(x',\omega^{\star})=-d_{g}(x,x')=\psi(x,x')\,.
$$
Similarly,
$$
\psi_{r}(x,\omega^{\star}+\pi)-\psi_{r}(x',\omega^{\star}+\pi)=d_{g}(x,x')=-\psi(x,x')\,.
$$ 
Therefore the values of the phase $\psi_{r}(x,\omega)-\psi_{r}(x',\omega)$
in the critical points are $\pm\psi(x,x')$. Note that this value is independent
of $r$. Coming back to (\ref{tos}) we apply the stationary phase formula to the
integral over $\omega$ and since $\pm\psi(x,x')$ is independent
of $r$, after the integration over $r$, we arrive at the claimed representation
for the kernel ${\mathcal K}(x,x')$. The non vanishing of $a^\pm (0,0)$ simply comes from the explicit representation we get from the stationary phase argument and the non negativeness of the function $a$ appearing in Theorem~\ref{sogge}. This reads
\begin{equation}
 a^{\pm}(0,x')\\
= \int_r \Psi(x,x',r) b_r(x,x', \omega^\star+ ( \pm \frac \pi 2 - \frac \pi 2)) dr
\end{equation}
where $\Psi(x,x',r)$ is a non vanishing function (the Jacobian determinant). 
  This completes the proof of Lemma~\ref{geom}.
\end{proof}
Let $\chi\in C_0^{\infty}(\R)$ such that 
\begin{equation}\label{zvezda}
\chi(x)=
\left\{
\begin{array}{l}
1,\quad |x|\leq 1/2,
\\
0,\quad |x|\geq 1
\end{array}
\right.
\end{equation}
and $\zeta\in C_{0}^{\infty}(\R)$ such that
$$
\zeta(x)=
\left\{
\begin{array}{l}
1,\quad |x|\leq 1/8,
\\
0,\quad |x|\geq 1/4\, .
\end{array}
\right.
$$
Then, since for $\lambda\geq 1$, 
$$
(1-\zeta)(\lambda x)(1-\chi)(\lambda^{1/3}x)=(1-\chi)(\lambda^{1/3}x),
$$ 
we can write
$$
1=\zeta(\lambda(t-\tau))+(1-\zeta)(\lambda(t-\tau))\,(1-\chi)(\lambda^{1/3}(t-\tau))+\chi(\lambda^{1/3}(t-\tau))\,.
$$
According to the above partition of the unity, we split the kernel as
\begin{equation}\label{eq.part}
K(t,\tau)=K_{1}(t,\tau)+K_{2}(t,\tau)+K_{3}(t,\tau)
\end{equation}
which in turn induces a natural splitting
of $T_{\lambda} T_{\lambda}^{\star}$ as
$$
T_{\lambda} T_{\lambda}^{\star} =L_{1}+L_{2}+L_{3}
$$
where $L_{j}$, $j=1,2,3$ denote the contributions of $K_{j}$ to 
$T_{\lambda} T_{\lambda}^{\star}$ respectively.
Using the Young inequality and a crude bound on $K(t,\tau)$, we obtain the
estimate
$$
\|L_{1}(f)\|_{L^{2}(\gamma)}\leq C\|\zeta(\lambda t)\|_{L^1(\mathbb{R})}
\|f\|_{L^2(\gamma)}
\leq C\lambda^{-1}\|f\|_{L^2(\gamma)}
$$
which is even better than the needed estimate. Next we estimate $L_2$ whose kernel in supported in the set 
$$ \{(t, \tau): \lambda^{-1} \lesssim |t- \tau| \lesssim \lambda^{-1/3}\}.$$
Using Lemma~\ref{geom}, we can write
$$
|K_{2}(t,\tau)|\leq
\frac{C\big|(1-\zeta)(\lambda(t-\tau))\,\chi(\lambda^{1/3}(t-\tau))\big|}
{(\lambda|t-\tau|)^{\frac{1}{2}}}\, .
$$
Using the Young inequality, we get
\begin{eqnarray*}
\|L_{2}(f)\|_{L^{2}(\gamma)}
& \leq &
\frac{C}{\lambda^{1/3}}
\Big\|
\frac{(1-\zeta)(\lambda t)\,\chi(\lambda^{1/3} t)}
{(\lambda^{1/3}|t|)^{\frac{1}{2}}}
\Big\|_{L^1(\mathbb{R})}\|f\|_{L^2(\gamma)}
\\
& \leq & C\lambda^{-\frac{2}{3}}\|f\|_{L^2(\gamma)}\, .
\end{eqnarray*}
It remains to estimate $L_3$ whose kernel is supported in the set
$$ \{(t, \tau): \lambda^{-1/3} \lesssim |t- \tau| \}.$$
For $t$, $\tau$ on the support of $(1-\chi)(\lambda^{1/3}(t-\tau))$, with the
notation of Lemma~\ref{geom}, we can write
$$
|b(x(t),x(\tau),\lambda^{-1}|x(t)-x(\tau)|^{-1})|\leq
C(\lambda|x(t)-x(\tau)|)^{-3/2}\leq C'\lambda^{-1}\, .
$$
Therefore, the contribution of $b$ to $L_3$ is bounded by $\lambda^{-1}\|f\|_{L^2( \gamma)}$ and we only need to estimate the contribution of the first term in the representation of
${\mathcal K}(x,x')$ given by Lemma~\ref{geom}.\par
Our next step is a general property of curves
passing through the origin of a geodesic coordinate system.
\begin{lemme}\label{geod}
Let 
$$
x(t)=(x_1(t),x_2(t)),\,\, t\sim 0  
$$ 
be a smooth curve in geodesic coordinate system with origin
$0$, parametrized by its arc length. Suppose also that $x(0)=0$ and $g(0)={\rm Id}$. Then
\begin{equation*}
|\dot{x}(0)|=1,\quad \langle\dot{x}(0),\ddot{x}(0)\rangle=0,\quad
|\ddot{x}(0)|^{2}=-\langle\dot{x}(0),\dddot{x}(0)\rangle,
\end{equation*}
where $\langle\cdot,\cdot\rangle$ is the $\R^2$ scalar product and $|\cdot|$
is the Euclidean norm.
\end{lemme}
\begin{proof}
We have that 
\begin{equation}\label{norm}
g(\dot{x}(t),\dot{x}(t))=1\, .
\end{equation}
and since $g(0)={\rm Id}$, we deduce that $|\dot{x}(0)|=1$.
Next, we differentiate (\ref{norm}) with respect to $t$ which gives
\begin{equation}\label{norm-bis}
2\sum_{i,j=1}^{2}g_{ij}(x(t))\ddot{x}_{i}(t)\dot{x}_{j}(t)+
\sum_{i,j,k=1}^{2}\frac{\partial g_{ij}}{\partial x_k}(x(t))\dot{x}_{i}(t)\dot{x}_{j}(t)\dot{x}_{k}(t)=0\, .
\end{equation}
Since we are working in geodesic coordinates the first derivatives of the metric are
vanishing at the origin.
Therefore by taking $t=0$ in (\ref{norm-bis}), we obtain that $\langle\dot{x}(0),\ddot{x}(0)\rangle=0$.
Since we are working in a geodesic coordinate system, we have the identity
\begin{equation}\label{base}
\sum_{j=1}^{2}g_{ij}(x(t))x_{j}(t)=\sum_{j=1}^{2}g_{ij}(0)x_{j}(t),\quad i=1,2\, .
\end{equation}
Differentiating three times (\ref{base}) with respect to $t$ and taking $t=0$
gives
$$
\sum_{j,k,s=1}^{2}
\frac{\partial^{2}g_{ij}}{\partial x_k \partial x_s}
(0)\,\dot{x}_{k}(0)\dot{x}_{j}(0)\dot{x}_{s}(0)=0\,
, \quad i=1,2.
$$
Therefore differentiating (\ref{norm-bis}) with respect to $t$ and taking $t=0$
gives
$$
|\ddot{x}(0)|^{2}=-\langle\dot{x}(0),\dddot{x}(0)\rangle\, .
$$
This completes the proof of Lemma~\ref{geod}\, .
\end{proof}
Next, we state an elementary lemma.
\begin{lemme}\label{taylor}
Let 
$$
x(t)=(x_1(t),x_2(t)),\,\, t\sim 0  
$$ 
be a smooth curve in $\R^2$ such that
\begin{equation}\label{degen}
x(0)=0,\quad |\dot{x}(0)|=1,\quad \langle\dot{x}(0),\ddot{x}(0)\rangle=0,\quad
|\ddot{x}(0)|^{2}=-\langle\dot{x}(0),\dddot{x}(0)\rangle,
\end{equation}
Then
\begin{equation}\label{nedelia}
|x(t)|^2= t^2\big(1-\frac{t^2}{12}\,|\ddot{x}(0)|^{2}+{\mathcal O}(|t|^3)\big)\,.
\end{equation}
\end{lemme}
\begin{proof}
Write the Taylor expansion
$$
x(t)=t\big(\dot{x}(0)+\frac{t}{2}\,\ddot{x}(0)+\frac{t^2}{6}\,\dddot{x}(0)+{\mathcal
  O}(|t|^3)\big)\, .
$$
We thus can write
\begin{eqnarray*}
|x(t)|^2= t^2
\big(
|\dot{x}(0)|^{2}+t\langle\dot{x}(0),\ddot{x}(0)\rangle+
\frac{t^2}{4}\,|\ddot{x}(0)|^{2}+\frac{t^2}{3}\,\langle\dot{x}(0),\dddot{x}(0)\rangle
+{\mathcal O}(|t|^3)\big)
\end{eqnarray*}
which, taking into account (\ref{degen}), clearly yields (\ref{nedelia}).
This completes the proof of Lemma~\ref{taylor}.
\end{proof}
With Lemma~\ref{taylor} in hand, we can now give a precise description of
$\psi(\gamma(t),\gamma(\tau))$. 
\begin{lemme}\label{phase}
For $\varepsilon\ll 1$, $\gamma(t),\gamma(\tau)\in U$, we can write
$$
\psi(\gamma(t),\gamma(\tau))=-|t-\tau|\Big(1-c(\tau)(t-\tau)^{2}+d(\tau,t-\tau)(t-\tau)^{3}\Big),
$$
where $c(\tau)$ and $d(\tau,t-\tau)$ are smooth functions such that 
$
c(\tau)\geq c>0\, .
$
\end{lemme}
\begin{proof}
Let us consider a geodesic coordinate system $x=(x_1,x_2)$ centered at $\gamma(\tau)$. In
this system the curve $\gamma$ is given by $\{x(s),\, s\sim 0 \}$ with
$x(0)=0$. 
Observe that 
$$
|x(t-\tau)|=-\psi(\gamma(t),\gamma(\tau))
$$
and therefore
we can apply Lemma~\ref{geod} and Lemma~\ref{taylor} to conclude that
\begin{equation}\label{back}
\psi^{2}(\gamma( t),\gamma(\tau))=|t-\tau|^{2}\Big(1-\frac{(t-\tau)^{2}}{12}\, |\ddot{x}(0)|^{2}+{\mathcal O}(|t-\tau|^3)\Big)\,.
\end{equation}
We have that in the considered coordinate system
$$
\gamma'\big|_{x(s)}=
\sum_{j=1}^{2}\dot{x}_j(s)\Big(\frac{\partial}{\partial x_j}\Big)_{x(s)}
$$
and by the definition
of the covariant derivative along $\gamma$ (see e.g. \cite{GaHuLa90}), 
we can write in the coordinates $(x_1,x_2)$, 
$$
\frac{D}{dt}\gamma'\Big|_{x(s)}
=\sum_{i=1}^{2}\left(
\ddot{x_{i}}(s)+\sum_{j,k=1}^{2}\Gamma_{jk}^{i}(x(s))\,\dot{x_j}(s)\dot{x_k}(s)
\right)\Big(\frac{\partial}{\partial x_i}\Big)_{x(s)},
$$
where $\Gamma_{jk}^{i}$ are the Christoffel symbols of the metric $g$.
Since we are working in a geodesic coordinate system, we have that
$
\Gamma_{jk}^{i}(0)=0\,.
$
Therefore 
$$
|\ddot{x}(0)|^{2}=g\Big(\frac{D}{dt}\gamma'(\tau)\, , \, \frac{D}{dt}\gamma^{'}(\tau)\Big)= 24 c(\tau)\neq 0
$$
thanks to the assumption (\ref{curvature}) on the curve $\gamma$. 
Consequently, coming back to (\ref{back}) and using that
$\psi$ is non-positive, we obtain the assertion of Lemma~\ref{phase}.
\end{proof}
In view Lemma~\ref{geom} and Lemma~\ref{phase}, it is now clear that
(\ref{1-bis}) is a consequence of the following oscillatory integral lemma.
\begin{lemme}\label{oscil}
Let $\chi\in C_0^{\infty}$ satisfying (\ref{zvezda}), $a\in C_0^{\infty}(\R^2)$
with support ${\mathcal O}(\varepsilon)$, $\varepsilon>0$. Let $\gamma\in
C^{\infty}(\R\times \R^{\star})$ be of the form
$$
\gamma(\tau,t-\tau)=\pm |t-\tau|\Big(1-c(\tau)(t-\tau)^{2}+d(\tau,t-\tau)(t-\tau)^{3}\Big),
$$
where $c(\tau)$ and $d(\tau,t-\tau)$ are smooth functions such that 
$
c(\tau)\geq c>0\, .
$
We define the linear map
$$
(K_{\lambda}f)(t)=\int_{-\infty}^{\infty}e^{i\lambda\gamma(\tau,t-\tau)}\frac{a(\tau,t-\tau)}{(\lambda|t-\tau|)^{\frac{1}{2}}}
\,\, (1-\chi)(\lambda^{1/3}(t-\tau))f(\tau)d\tau\,.
$$
Then for $\varepsilon\ll 1$,
$$
\|K_{\lambda}f\|_{L^2}\leq C\lambda^{-\frac{2}{3}}\|f\|_{L^2}\,\, .
$$
\end{lemme} 
\begin{proof}
As a startup and to enlighten  the exposition below, let us first indicate how to prove the result in the model case 

\begin{equation}
\gamma( \tau, t- \tau)= (t- \tau) + c (t-\tau)^3,\qquad
a(\tau, t- \tau) =1
\end{equation}
In that case our operator is (modulo a conjugation by $e^{i\lambda s}$ which is an isometry on $L^2$) a convolution operator with convolution kernel 
$$ k(r) = e^{i \lambda r^3} \frac{(1-\chi)(\lambda^{1/3} r)} { \lambda |r|^{1/2}}$$
and consequently the norm of the operator is equal to 
$$ \| \widehat{k}\|_{L^\infty}= \lambda^{-2/3} \sup_\rho |\widehat{\widetilde{k}}|(\lambda^{-1/3} \rho)= \lambda^{-2/3} \sup_\rho |\widehat{\widetilde {k}}|(\rho) $$
with 
$$ \widetilde{k}(r) = e^{i r^3}  \frac{(1-\chi)( r)} { |r|^{1/2}}$$
having a bounded Fourier transform, as can be seen by a straightforward application of van der Corput lemma.\par
 We now come back to the general case. 
We shall only consider the case when the sign in front of $|t-\tau|$ in the
definition of $\gamma(\tau,t-\tau)$ is minus, the case of sign plus being
similar. We split the kernel of $K_\lambda$ into two parts supported respectively in the sets $\{t\geq \tau\}$ and $\{t \leq \tau\}$ which yields 
$$
K_{\lambda}=K_{\lambda}^{+}+K_{\lambda}^{-}.
$$
With 
$$ \gamma_{1}(\tau,t-\tau)=\gamma(\tau,t-\tau)+|t-\tau|= c(\tau)( t- \tau)^3 - d(\tau, t- \tau) ( t- \tau)^4, $$
we can write
$$
(K_{\lambda}^{+}f)(t)=e^{-i\lambda t}
\int_{t\geq\tau}\,
e^{i\lambda\gamma_{1}(\tau,t-\tau)}\frac{a(\tau,t-\tau)}{(\lambda|t-\tau|)^{\frac{1}{2}}}
\,\,(1-\chi)(\lambda^{1/3}(t-\tau))e^{i\tau \lambda}f(\tau)d\tau\,.
$$
It is now clear that  it suffices to establish the estimate
$$
\|\widetilde{K}_{\lambda}^{+}f\|_{L^2}\leq C\lambda^{-\frac{2}{3}}\|f\|_{L^2}\,\, ,
$$
where
$$
(\widetilde{K}_{\lambda}^{+}f)(t)=
\int_{t\geq\tau}\,
e^{i\lambda\gamma_{1}(\tau,t-\tau)}\frac{a(\tau,t-\tau)}{(\lambda|t-\tau|)^{\frac{1}{2}}}
\,\,(1-\chi)(\lambda^{1/3}(t-\tau))f(\tau)d\tau\,.
$$
Let us compute the Fourier transform
$$
\widehat{\widetilde{K}_{\lambda}^{+}f}(\xi)=
\int_{-\infty}^{\infty}\int_{t\geq\tau}
\exp\Big(i\lambda\gamma_{1}(\tau,t-\tau)-it\xi\Big)\frac{a(\tau,t-\tau)}{(\lambda|t-\tau|)^{\frac{1}{2}}}
(1-\chi)(\lambda^{1/3}(t-\tau))f(\tau)d\tau dt\,\, .
$$
We now exchange the order of integration and then we perform a change of
variables
$(\tau,t)\rightarrow (\tau,z)$
where $t=\tau+z/\lambda^{1/3}$. Therefore
\begin{multline}
\widehat{\widetilde{K}_{\lambda}^{+}f}(\xi)\\
=\lambda^{-\frac{2}{3}}
\int_{-\infty}^{\infty}\int_{1/2}^{\infty}
\exp\Big(i\lambda\gamma_{1}\big(\tau,\frac{z}{\lambda^{1/3}}\big)-i\tau\xi-i\frac{z\xi}{\lambda^{1/3}}\Big)
a\Big(\tau,\frac{z}{\lambda^{1/3}}\Big)\widetilde{\chi}(z)f(\tau)dz d\tau,
\end{multline}
where 
$$
\widetilde{\chi}(z)=\frac{(1-\chi)(z)}{z^{1/2}}\, .
$$
We can therefore write
$$
\widehat{\widetilde{K}_{\lambda}^{+}f}(\xi)=\lambda^{-\frac{2}{3}}\,{\mathcal K}f(\xi)
$$
where
$$
{\mathcal K}f(\xi)=\int_{-\infty}^{\infty}e^{-i\tau\xi}A_{\lambda}(\tau,\xi)f(\tau)d\tau
$$
with
\begin{equation}\label{eq.4.16b}
A_{\lambda}(\tau,\xi)=\int_{1/2}^{\infty}
\exp\Big(i\big[\lambda\gamma_{1}\big(\tau,\frac{z}{\lambda^{1/3}}\big)-\frac{\xi z}{\lambda^{1/3}}\big]\Big)
a\Big(\tau,\frac{z}{\lambda^{1/3}}\Big)\widetilde{\chi}(z)dz\,.
\end{equation}
We therefore reduced the matters to establishing a uniform $L^2$ bound on the
linear map ${\mathcal K}$. Let us denote by $\Phi_{\lambda}(\tau,\xi,z)$ the
phase function in the definition of $A_{\lambda}(\tau,\xi)$, i.e.
$$
\Phi_{\lambda}(\tau,\xi,z)=
\lambda
\gamma_{1}\big(\tau,\frac{z}{\lambda^{1/3}}\big)-\frac{\xi z}{\lambda^{1/3}} \, .
$$
Coming back to the definition of $\gamma_{1}$, we can write
$$
\Phi_{\lambda}(\tau,\xi,z)=c(\tau)z^{3}-\xi\, \frac{z}{\lambda^{1/3}}
-d\big(\tau,\frac{z}{\lambda^{1/3}}\big)\,\frac{z^{4}}{\lambda^{1/3}}\, .
$$
Hence
$$
\frac{\partial \Phi_{\lambda}}{\partial z}(\tau,\xi,z)=
3c(\tau)z^{2}-
\frac{\xi}{\lambda^{1/3}}
-d\big(\tau,\frac{z}{\lambda^{1/3}}\big)
\frac{z}{\lambda^{1/3}}\,4z^2-
\frac{\partial d}{\partial z}\big(\tau,\frac{z}{\lambda^{1/3}}\big)\Big(\frac{z}{\lambda^{1/3}}\Big)^{2}z^{2}\,\, .
$$
Therefore for $\varepsilon\ll 1$ and $(z, \tau)$ on the support of
$a\big(\tau,\frac{z}{\lambda^{1/3}}\big)$
\begin{equation}\label{dev1}
\Big|
\frac{\partial \Phi_{\lambda}}{\partial z}(\tau,\xi,z)
\Big|\geq cz^{2}\, ,
\end{equation}
provided $\xi\leq 0$. 
On the other hand for every $\xi\in\R$, every $k\geq 2$, every $(z,\tau)$
on the support of $a\big(\tau,\frac{z}{\lambda^{1/3}}\big)$, we have the bound
\begin{equation}\label{eq.estphase}
\Big|
\frac{\partial^{k} \Phi_{\lambda}}{\partial z^{k}}(\tau,\xi,z)
\Big|\leq Cz^{\max(3-k,0)} \, .
\end{equation}
Therefore, we can integrate by parts with the operator 
$$ \frac{ 1 } { i \frac{ \partial \Phi_\lambda} { \partial z}} { \frac {\partial} { \partial z}}.$$
Each such integration by parts gains $\frac 1 {c z^3}$. Since derivatives with respect to $\tau$ give only (fixed) powers of $z$, we obtain that for $\xi\leq 0$,
$k=0,1,\dots$,
\begin{equation}\label{dev2}
|\partial_{\tau}^{k}A_{\lambda}(\tau,\xi)|\leq C_{k}\,\, .
\end{equation}
Notice that if $|\xi|\leq C\lambda^{1/3}$ then we can put the slowly oscillating factor
$e^{-i\xi z/\lambda^{1/3}}$ in the amplitude $\widetilde{\chi}(z)$, and, the
new phase satisfies the bound (\ref{dev1})  and similar integration by parts gaining now $\frac 1 { c z^2}$, we can still 
achieve the bound (\ref{dev2}). 
\par
We obtain that for every $\Lambda>0$, every $k=0,1,\dots$ there
exists $C>0$ such that for every $\xi\leq \Lambda \lambda^{1/3}$,
\begin{equation}\label{vajna}
|\partial_{\tau}^{k}A_{\lambda}(\tau,\xi)|\leq C\, .
\end{equation}
Let $\chi_{0}\in C_0^{\infty}(\R)$ be such that
$$
\chi_0(x)=
\left\{
\begin{array}{l}
0,\quad |x|\leq 1/2,
\\
1,\quad |x|\geq 1\, .
\end{array}
\right.
$$
For $\Lambda>0$ to be fixed later, we consider the splitting
$$
{\mathcal K}={\mathcal K}_{1}+{\mathcal K}_{2},
$$
where
\begin{equation*}
({\mathcal K}_{1}f)(\xi)  =  (1-\chi_0)\Big(\frac{\xi}{\Lambda \lambda^{1/3}}\Big)
({\mathcal K}f)(\xi)
\end{equation*}
Observe that 
$$({\mathcal K}_1f)(\xi)= \int_{-\infty}^{\infty} (1-\chi_0)\Big(\frac{\xi}{\Lambda \lambda^{1/3}}\Big)
\widehat{A_{\lambda}}\big(\xi-\sigma,\xi\big)\hat{f}(\sigma)\frac{d\sigma}{2\pi}\,\, ,
$$
where $\widehat{A_{\lambda}}$ denotes the Fourier transform of $A_{\lambda}$
with respect to the first variable. Thanks to (\ref{vajna}), and using the
compactness of the support of $A_{\lambda}(\tau,\xi)$ with respect to $\tau$,
we obtain that for every $N\in\N$ there exist $C_N$ such that
\begin{equation}\label{Young}
| (1-\chi_0)\Big(\frac{\xi}{\Lambda \lambda^{1/3}}\Big)\widehat{A_{\lambda}}\big(\xi-\sigma,\xi\big)|\leq
\frac{C_N}{(1+|\xi-\sigma|)^{N}}\,\, .
\end{equation}
Estimate (\ref{Young}) with $N=2$ and the Young inequality
imply the uniform (with respect to $\lambda$) boundedness of 
${\mathcal K}_{1}$ on $L^2(\R)$. 
\par
It remains to estimate ${\mathcal K}_{2}$ as an operator on $L^{2}(\R)$.
For $\xi>0$, we consider a change of variable
$$
z\mapsto \Big(\frac{\xi}{\lambda^{1/3}}\Big)^{1/2}z
$$
which allows to write
\begin{equation}\label{int}
A_{\lambda}(\tau,\xi)=\Big(\frac{\xi}{\lambda^{1/3}}\Big)^{\frac{1}{2}}
\int
\exp\Big(i\big(\frac{\xi}{\lambda^{1/3}}\big)^{\frac{3}{2}}\phi_{\lambda}(\tau,\xi,z)\Big)
a\Big(\tau,\frac{\sqrt{\xi}\,z}{\sqrt{\lambda}}\Big)
\widetilde{\chi}\big(\frac{\sqrt{\xi}\,z}{\lambda^{1/6}}\big)dz.
\end{equation}
Remark that due to support considerations, in the integral above 
$$
z\geq \frac{1}{2}\left(\frac{\xi}{\lambda^{1/3}}\right)^{-\frac{1}{2}}
$$ 
and
\begin{equation}\label{phi}
\phi_{\lambda}(\tau,\xi,z)=c(\tau)z^{3}-z-
d\big(\tau,\frac{\sqrt{\xi}\,z}{\sqrt{\lambda}}\big)\,\frac{\sqrt{\xi}\,z}{\sqrt{\lambda}}\,z^{3}\,\, .
\end{equation}
Now we set 
$$\omega= \Bigl( \frac{\xi} { \lambda^{1/3}}\Bigr) ^{3/2}$$ and we intend to apply the stationary phase formula ($\omega$ being the large parameter) to the integral~\eqref{int}. Consider the equation
\begin{equation}\label{equation}
3c(\tau)z^{2}-1+{\mathcal O}(\varepsilon z^{2})=0
\end{equation}
satisfied by the critical points with respect to $z$ of the phase
$\phi_{\lambda}(\tau,\xi,z)$ on the support of $a\Big(\tau,\frac{\sqrt{\xi}\,z}{\sqrt{\lambda}}\Big)$.
For $\varepsilon\ll 1$, the equation
(\ref{equation}) has a unique positive root $z=z_{\lambda}(\tau,\xi)$ which
satisfies
\begin{equation}\label{ineq}
0<c_1\leq z_{\lambda}(\tau,\xi)\leq c_{2}\,\, .
\end{equation}
In addition, the critical point
$z_{\lambda}(\tau,\xi)$ is non degenerate
$$ \exists c>0; c\leq \frac{\partial^2 \phi_\lambda} { \partial z^2} (z_{\lambda}(\tau,\xi)) \leq \frac 1 c.$$
Since
\begin{equation}\label{i}
\frac{\partial\phi_{\lambda}}{\partial z}(\tau,\xi,z_{\lambda}(\tau,\xi))=0\, .
\end{equation}
Differentiating (\ref{i}) with respect to $\tau$ gives
\begin{equation}\label{ii}
\frac{\partial^{2}\phi_{\lambda}}{\partial z\partial\tau}(\tau,\xi,z_{\lambda}(\tau,\xi))
+
\frac{\partial z_{\lambda}}{\partial \tau}(\tau,\xi)\,
\frac{\partial^{2}\phi_{\lambda}}{\partial z^{2}}(\tau,\xi,z_{\lambda}(\tau,\xi))
=0\, .
\end{equation}
Coming back to (\ref{phi}) and using~\eqref{ineq}, we obtain that
$$
\left|
\frac{\partial z_{\lambda}}{\partial \tau}(\tau,\xi)
\right|
\leq C\sqrt{\frac{\xi}{\lambda}}
$$

Similarly, differentiating (\ref{ii}) with respect to $\tau$, we get for
$k\geq 1$,
\begin{equation}\label{iii}
\left|
\frac{\partial^{k} z_{\lambda}}{\partial \tau^{k}}(\tau,\xi)
\right|\leq C_{k}\, .
\end{equation}
 Next, we notice that
\begin{equation}\label{eq.estsymb}
\left|
\big(\frac{\partial}{\partial \tau}  \big)^{p}
\big(\frac{\partial}{\partial z}  \big)^{q}
a\Big(\tau,\frac{\sqrt{\xi}\,z}{\sqrt{\lambda}}\Big)
\widetilde{\chi}\big(\frac{\sqrt{\xi}\,z}{\lambda^{1/6}}\big)
\right|
\leq C_{p,q}\Big(\frac{\xi}{\lambda^{1/3}}\Big)^{\frac{q}{2}}\,\, .
\end{equation}
But for $\omega \gg 1$ we have that
$$
\Big(\frac{\xi}{\lambda^{1/3}}\Big)^{{q}}\ll
\omega^q 
$$
Therefore, for $\xi/\lambda\leq C$, we can apply the stationary phase formula to sufficiently high order (we recall that in the expansion given by the stationary phase formula the $k$-th term in $\omega^{-k}$ involves derivatives of order $2k$ of the symbol, see~\cite[Th. 7.7.1]{Ho83}), and we get
$$
A_{\lambda}(\tau,\xi)=
\exp\Big(
i\omega
V_{\lambda}(\tau,\xi)
\Big)q_{\lambda}(\tau,\xi)
$$
where $q_{\lambda}(\tau,\xi)$ satisfies
$$
|\partial_{\tau}^{k}q_{\lambda}(\tau,\xi)|\leq C_{k}
\Big(\frac{\xi}{\lambda^{1/3}}\Big)^{-\frac{1}{4}} \, , \quad k=0,1,2,3
$$
and $V_\lambda(\tau, \xi)$ is the value of the phase $\phi_\lambda( \tau, \xi, z)$ at the critical point $z= z_\lambda( \tau, \xi)$, hence satisfies
\begin{equation}\label{V}
|\partial_{\tau}^{k}V_{\lambda}(\tau,\xi)|\leq C_{k}\, ,\quad k=0,1,2,3.
\end{equation}
Moreover, if $\frac \xi \lambda \gg 1$, then the critical point $z_\lambda(\tau, \xi)$ lies far from the support of the function $a$ and consequently, by integrating by parts and using~\eqref{eq.estphase}~\eqref{eq.estsymb}, we obtain a contribution rapidly decaying with respect to $\lambda$.
Observe that the condition $\omega \gg 1$ can be achieved by taking
$\Lambda\gg 1$
in the definition of the cut-off $\chi_{0}$.\par
Let $\chi_{1}\in C_{0}^{\infty}(\R)$ be equal to one in a  neighborhood of
zero. Consider the splitting
$$
{\mathcal K}_{2}={\mathcal K}_{21}+{\mathcal K}_{22},
$$
where
$$
({\mathcal K}_{21}f)(\xi)=
(1-\chi_1)\Big(\frac{\xi}{\delta \lambda}\Big)
({\mathcal K}_{2}f)(\xi)\, .
$$
Notice that for $\xi$ on the support of
$(1-\chi_1)\Big(\frac{\xi}{\delta \lambda}\Big)$ one has
\begin{equation}\label{raz}
|\xi|\gtrsim \delta \lambda\,\, .
\end{equation}
On the other hand, on the domain of integration in (\ref{int}),
\begin{equation}\label{dva}
\left|
\frac{\sqrt{\xi}\, z}{\sqrt{\lambda}}
\right|\lesssim \varepsilon\, .
\end{equation}
Therefore, if $\delta= \kappa \varepsilon^2$ with $\kappa $ large enough, we obtain that (\ref{raz}) and (\ref{dva}) imply $|z|\leq c_{1}/2$.
Hence, in view of (\ref{ineq}), for $\xi$ on the support of
$(1-\chi_1)\Big(\frac{\xi}{\delta\lambda}\Big)$
the phase function $\phi_{\lambda}(\tau,\xi,z)$ has no critical point in $z$
and we can evaluate ${\mathcal K}_{21}$ as we did for ${\mathcal K}_{1}$ thanks to the rapid decay (in $\lambda$) satisfied by
${A}_{\lambda}(\tau,\xi)$.\par
It remains to deal with ${\mathcal K}_{22}$.
For that purpose, we make another appeal to a duality argument.
We need to establish the uniform (with respect to $\lambda$) boundedness on
$L^2(\R)$ of the operator
$$
({\mathcal K}_{22}f)(\xi)=
\chi_{0}\Big(\frac{\xi}{\Lambda\lambda^{1/3}}\Big)
\chi_1\Big(\frac{\xi}{\delta \lambda}\Big)
\int_{-\infty}^{\infty}
e^{-i\tau\xi}
e^{i\omega\,V_{\lambda}(\tau,\xi)}
q_{\lambda}(\tau,\xi)f(\tau)d\tau\, .
$$
By duality one therefore needs to study the uniform $L^2$ boundedness of
$
{\mathcal K}_{22}{\mathcal K}_{22}^{\star},
$
where ${\mathcal K}_{22}^{\star}$ is the formal adjoint of ${\mathcal K}_{22}$
with respect to the $L^2(\R)$ inner product. Write
$$
({\mathcal K}_{22}{\mathcal K}_{22}^{\star}f)(\xi)=
\int_{-\infty}^{\infty}S(\xi,\sigma)f(\sigma)d\sigma,
$$
where
\begin{multline*}
S(\xi,\sigma)=
\chi_{0}\Big(\frac{\xi}{\Lambda\lambda^{1/3}}\Big)
\chi_1\Big(\frac{\xi}{\delta\lambda}\Big)
\chi_{0}\Big(\frac{\sigma}{\Lambda\lambda^{1/3}}\Big)
\chi_1\Big(\frac{\sigma}{\delta \lambda}\Big)\times
\\
\times
\int_{-\infty}^{\infty}
e^{-i\tau(\xi-\sigma)}
\exp\Big(
i\big[\omega V_{\lambda}(\tau,\xi)
-
\frac{\sigma^{3/2}}{\lambda^{1/2}}\,V_{\lambda}(\tau,\sigma)
\big]
\Big)
q_{\lambda}(\tau,\xi)\overline{q_{\lambda}}(\tau,\sigma)d\tau\, .
\end{multline*}
We next write
\begin{equation}\label{Tay}
\omega V_{\lambda}(\tau,\xi)
-
\frac{\sigma^{3/2}}{\lambda^{1/2}}\,V_{\lambda}(\tau,\sigma)
=\frac{\xi^{3/2}}{\lambda^{1/2}}\,V_{\lambda}(\tau,\xi)
-
\frac{\sigma^{3/2}}{\lambda^{1/2}}\,V_{\lambda}(\tau,\sigma)
=
(\xi-\sigma)R_{\lambda}(\tau,\xi,\sigma)\, .
\end{equation}
The relevant fact about $R_{\lambda}$ is that for $\varepsilon\ll 1$, the map
$$
\tau\mapsto \tau-R_{\lambda}(\tau,\xi,\sigma)
$$
is a small perturbation of the identity and thus the corresponding integral
operator is $L^2$ bounded. Let us now give the precise argument.
$$
V_{\lambda}(\tau,\xi)=\phi_{\lambda}(\tau,\xi,z_{\lambda}(\tau,\xi))
$$
and therefore, by invoking (\ref{i}), we obtain that
$$
\frac{\partial V_{\lambda}}{\partial \xi}(\tau,\xi)
=
\frac{\partial\phi_{\lambda}}{\partial \xi}(\tau,\xi,z_{\lambda}(\tau,\xi))\,.
$$
Coming back to the definition of $\phi_{\lambda}$, using (\ref{ineq}) and
(\ref{iii}), we obtain that
\begin{equation}\label{v}
\left|
\Big(
\frac{\partial}{\partial \tau}
\Big)^{k}\,
\frac{\partial V_{\lambda}}{\partial \xi}(\tau,\xi)
\right|
\leq C_{k}\,\xi^{-\frac{1}{2}}\,\lambda^{-\frac{1}{2}}\,.
\end{equation}
Since $\delta \lesssim \epsilon^2$, estimates (\ref{V}) and (\ref{v}) imply the bound
$$
\left|
\Big(
\frac{\partial}{\partial \tau}
\Big)^{k}\,
\Big(
\frac{\partial}{\partial \xi}
\Big)
\Big[
\frac{\xi^{3/2}}{\lambda^{1/2}}\,
V_{\lambda}(\tau,\xi)
\Big]
\right|
\leq C\Big(\sqrt{\frac{\xi}{\lambda}}+\frac{\xi}{\lambda}\Big)
\leq C\varepsilon,\quad k=0,1,2,3
$$
for $\xi$ on the support of 
$\chi_1\Big(\frac{\xi}{\delta\lambda}\Big)$.
Coming back to (\ref{Tay}), we obtain
$$
|\partial_{\tau}^{k}\,R_{\lambda}(\tau,\xi,\tau)|\leq C\varepsilon\,
,\quad k=0,1,2,3.
$$
For $\varepsilon\ll 1$, after two integrations by parts in $\tau$, we obtain the bound
$$
|S(\xi,\sigma)|\leq \frac{C}{(1+|\xi-\sigma|)^{2}}
$$
and therefore using the Young inequality we obtain that
$
{\mathcal K}_{22}{\mathcal K}_{22}^{\star}
$
is bounded on $L^2(\R)$, uniformly in $\lambda$.
This completes the proof of Lemma~\ref{oscil} and consequently of~(\ref{1-bis}).
\end{proof}  
\section{On the optimality of (\ref{1}) and (\ref{1-bis})}\label{sec.opt}
\subsection{Optimality for the spectral projectors on any manifold}
In this section we prove lower bounds for the smoothed spectral projectors $\chi_\lambda$ involved in Theorem~3. In section~\ref{sec.5.2} we will use these bounds to get the optimality of~\eqref{1} and ~\eqref{1-bis} in the particular case of the standard sphere. Let us first study the case $p\geq4$.
\begin{lemme} Let $p\geq 4$. For any smooth curve $\gamma$, there exists $c>0$ such that for any $\lambda \geq 1$,
\begin{equation}\label{eq.opt1}
\| \chi_\lambda \|_{\mathcal{L}(L^2(M); L^p( \gamma)))}\geq c \lambda ^{\frac 1 2 - \frac 1 p} 
\end{equation}
\end{lemme}
\begin{proof}
Indeed since 
 $$\| \chi_\lambda \|^2_{\mathcal{L}(L^2(M); L^p( \gamma)))}=\| \chi_\lambda\chi_\lambda^* \|_{\mathcal{L}(L^{p'}(\gamma); L^p( \gamma)))}$$
and since 
$$\| R_\lambda \|_{\mathcal{L}(L^2(M); L^p( \gamma)))}\leq C\lambda^{-1/2},$$
it suffices to prove
$$\|T_{\lambda}T_\lambda^\star \|_{\mathcal{L}(L^{p'}(\gamma); L^p(\gamma))} \geq c \lambda ^{-\frac 2 p}.$$ 
But if $|t-s|\leq 2\epsilon\lambda^{-1}$ then the oscillatory factor $e^{i\lambda[\psi(x(t),y)- \psi(x(s),y)]}$ does not oscillate any more and if $\epsilon>0$ is chosen small enough, the kernel ${K}$ of the operator $T_{\lambda}T_\lambda^\star$ satisfies
\begin{multline*}\Re(K(t,s))=\Re(\mathcal{K}(x(t),x(s)))\\
= \int_V \Re (e^{i\lambda[\psi(x(t),y)- \psi(x(s),y)]}) a(x(t),y) a(x(s),y) dy \geq \alpha >0.
\end{multline*}
As a consequence, choosing as test functions $f(s)=\lambda ^{\frac 1 {p'}} \phi\bigl( \frac {\lambda s} \epsilon\bigr)$ (with $\phi\in C^\infty_0 (-1,1)$ real and non negative), we obtain (for a constant $c>0$)
\begin{equation*}
\begin{aligned}
\Re \Bigl(1_{|t|\leq \epsilon\lambda^{-1}}T_{\lambda}T_\lambda^\star f\Bigr) &\geq c1_{|t|\leq \epsilon\lambda^{-1}}\int  \lambda ^{\frac 1 {p'}} \phi\Bigl( \frac {\lambda s} \epsilon \Bigr)ds\\
&\geq c1_{|t|\leq \epsilon\lambda^{-1}}  \lambda ^{\frac 1 {p'}-1} 
\end{aligned}
\end{equation*}
and 
$$\|T_{\lambda}T_\lambda^\star f\|_{L^p} \geq c \lambda ^{\frac 1 {p'}-1- \frac 1 p}=c \lambda ^{-\frac 2 p}$$
\end{proof}
We now turn to the case $2\leq p \leq 4$ {\em if $\gamma$ is a geodesic}.
\begin{lemme} Let $2\leq p\leq 4$. For any geodesic $\gamma$, there exist $c>0$ such that for any $\lambda \geq 1$,
\begin{equation}\label{eq.opt2}
\| \chi_\lambda \|_{\mathcal{L}(L^2(M); L^p( \gamma)))}\geq c \lambda ^{\frac 1 4} 
\end{equation}
\end{lemme}
\begin{proof}
As in the previous case, this lemma is equivalent to
$$\|T_{\lambda}T_\lambda^\star \|_{\mathcal{L}(L^{p'}(\gamma); L^p(\gamma))} \geq c \lambda ^{-\frac 1 2}.$$
Assume that the geodesic is parametrized by arc length. The critical value of the phase in the expression of the kernel of the operator $T_\lambda T_\lambda^\star$ is
$$ \psi(x(t), x(s))= |t-s|$$
and coming back to the expression of the kernel given by Lemma~\ref{geom} we have for $C\lambda^{-1} \leq |t-s|$,
\begin{multline}
K(t,s) = \mathcal{K}(x(t), x(s))\\
 = e^{i\lambda|t-s|}\frac {a^+(t,s)}{(\lambda|t-s|)^{1/2}}+ e^{-i\lambda|t-s|}\frac {a^-(t,s)}{(\lambda|t-s|)^{1/2}} + \mathcal{O} \Bigl( \frac 1 {(\lambda|t-s|)^{3/2}}\Bigr)
\end{multline}
where $a^\pm$ is non negative and does not vanish for $t,s$ in a neighborhood of $(0,0)$. 
As a consequence, choosing as test functions $f(s)=e^{i\lambda s} \phi(  \frac {s- 2 \epsilon} \epsilon)$ (with $\phi\in C^\infty_0 (-1,1)$ real and non negative), we obtain
\begin{multline}\label{eq.geodopt}
\Re \Bigl(e^{-i\lambda t}1_{4 \epsilon\leq  t\leq 5\epsilon }T_{\lambda}T_\lambda^\star f\Bigr) \\
\geq 1_{ 4 \epsilon\leq  t\leq 5\epsilon }\Re \Bigl(\int e^{-i\lambda t} K(t,s) e^{i\lambda s} \phi( \frac  {s- 2 \epsilon} \epsilon)ds\Bigr)\\
\geq  1_{ 4 \epsilon\leq  t\leq 5\epsilon } \int_{\epsilon \leq s \leq 2 \epsilon} \frac{a^+(t,s)}{(\lambda(t-s))^{1/2}}\phi( \frac {s- 2 \epsilon} \epsilon)  ds\\
 - \Bigl|\int_{\epsilon \leq s \leq 2 \epsilon} e^{-2i\lambda (t-s)}\frac{a^-(t,s)}{(\lambda(t-s))^{1/2}} \phi( \frac {s- 2 \epsilon} \epsilon) ds\Bigr| - \mathcal{O}(\lambda^{-3/2})
\end{multline}
We can integrate by parts with respect to $s$ in the second integral in the right hand side of~\eqref{eq.geodopt} and gain any power of $\lambda^{-1}$. As a consequence, the main contribution is the first term and we obtain
$$\|T_{\lambda}T_\lambda^\star f\|_{L^p} \geq \frac c{\lambda ^{1/2}} \qquad(c  >0),$$
which completes the proof. \end{proof}

We finally turn to the case $2\leq p \leq 4$ {\em if $\gamma$ is geodesically curved}.
\begin{lemme}\label{lem.optbis} Let $2\leq p\leq 4$. For any geodically curved $\gamma$, there exists $c>0$ such that for any $\lambda \geq 1$,
\begin{equation}\label{eq.opt3}
\| \chi_\lambda \|_{\mathcal{L}(L^2(M); L^p( \gamma)))}\geq c \lambda ^{{\frac 1 3 - \frac 1 {3p}} } 
\end{equation}
\end{lemme}
\begin{proof}
As before this lemma is equivalent to
$$\|T_{\lambda}T_\lambda^\star \|_{\mathcal{L}(L^{p'}(\gamma); L^p(\gamma))} \geq c \lambda ^{-\frac 1 3 -\frac 2 {3p}}.$$
In that case, we use the expression of the phase obtained in~\eqref{back} and  choose as test functions $f(s)= e^{i\lambda s}\lambda ^{\frac 1 {3p'}} \phi( \frac {\lambda^{1/3} s-2\epsilon} \epsilon)$ (with $\phi\in C^\infty_0 (-1,1)$ real and non negative) and obtain 
\begin{multline}
\Re \Bigl(e^{-i\lambda t}1_{4\epsilon \lambda^{-1/3}\leq t\leq 5\epsilon\lambda^{-1/3}}T_{\lambda}T_\lambda^\star f \Bigr)\\
\geq 1_{4\epsilon \lambda^{-1/3}\leq t\leq 5\epsilon\lambda^{-1/3}}\Re \Bigl(\int e^{-i\lambda t}K(t,s) e^{i \lambda s} \lambda ^{\frac 1 {3p'}} \phi \Bigl( \frac {\lambda^{1/3} s-2\epsilon} \epsilon\Bigr)ds\Bigr) 
\end{multline}
But on the support in $s$ of the function $f$ and for $4\epsilon \lambda^{-1/3}\leq t\leq 5\epsilon\lambda^{-1/3}$, we have 
$$\epsilon \lambda^{-1/3} \leq s \leq 3\epsilon \lambda^{-1/3} \leq 4\epsilon \lambda^{-1/3} \leq t \leq 5 \epsilon \lambda^{-1/3} $$ and consequently we can use the description of the kernel given in Lemma~\ref{geom} and  according to Lemma~\ref{phase}, for $\epsilon>0$ small enough
$$ \Re ( e^{-i\lambda t} e^{-i \lambda \psi(t,s) } e^{i\lambda s})=\Re \Bigl(e^{-ic(s) \lambda  (t-s) ^3 (1+ \mathcal{O}( \lambda^{-1/3}))}\Bigr) \geq \frac 1 2 $$
We deduce, since $a^+(0,0)$ is positive,
 \begin{multline}
\Re \Bigl(e^{-i\lambda t}1_{4\epsilon \lambda^{-1/3}\leq t\leq 5\epsilon\lambda^{-1/3}}T_{\lambda}T_\lambda^\star f \Bigr)\\
\begin{aligned}
&\geq c 1_{ 4\epsilon\lambda^{-1/3}\leq t\leq 5\epsilon\lambda^{-1/3}}\int  \phi \Bigl( \frac {\lambda^{1/3} s-2\epsilon} \epsilon\Bigr)\times \frac 1 { (\lambda (t-s))^{1/2}}ds\\
& \qquad \hfill -\Bigl| \int  e^{-i\lambda ((t-s)+ \psi(t,s))} \frac{a^-(t,s)}{(\lambda(t-s))^{1/2}} \phi \Bigl( \frac {\lambda^{1/3} s- 2 \epsilon} \epsilon\Bigr)ds\Bigr| \\
&\qquad \hfill - \mathcal{O}( \lambda^{-1}) \int  \phi \Bigl( \frac {\lambda^{1/3} s-2\epsilon} \epsilon\Bigr)ds\\
&\geq c\lambda^{-1/3}1_{ 4\epsilon\lambda^{-1/3}\leq t\leq 5\epsilon\lambda^{-1/3}}\int   \phi \Bigl( \frac {\lambda^{1/3} s-2\epsilon} \epsilon\Bigr)ds- \mathcal{O}( \lambda^{\frac1 {3p'}-2})
\end{aligned}\end{multline}
where as in the previous case we integrated by parts with respect to the variable $s$ to gain powers of $\lambda^{-1}$ in the integral involving $a^-$ (each such integration by parts gains  $\lambda^{-1}$ and looses $\lambda^{1/3}$ due to the powers of $(t-s)$). As a consequence
$$\|T_{\lambda}T_\lambda^\star f\|_{L^p} \geq c  \lambda ^{-\frac 2 3 - \frac 1 {3p}}=\lambda ^{\frac 1 {3p'}-\frac 2 3 - \frac 1 {3p}}\|f\|_{L^{p'}}=c  \lambda ^{-\frac 1 3-\frac 2 {3p}}\|f\|_{L^{p'}}$$
which completes the proof. \end{proof}
\subsection{From spectral projectors to exact eigenfunctions on spheres}\label{sec.5.2}
\begin{proof}We first remark that our lower bound for the norm of the smoothed spectral projector $\chi_\lambda$ from $L^2(M)$ to $L^p(\gamma)$ gives a lower bound for the norm of the spectral projector
$$ \Pi_\lambda= 1_{\sqrt{- {\mathbf{\Delta}}} -\lambda \in [0, \frac 1 2)}$$ in the following sense:
\begin{equation}\label{eq.seq} \exists c>0, \exists \lambda_n\rightarrow+ \infty , f_n \in L^2(M); \|\Pi_{\lambda_n} f_n \|_{L^p( \gamma)} > c \lambda_n^{\check{\delta}(p)}\|f_n\|_{L^2(M)}
\end{equation}
where $\check{\delta}(p) = \widetilde{\delta}(p)$ or $\check{\delta}(p) = \delta(p)$ according whether the curved is geodesically curved or not. 
Indeed, if the converse were true,
$$  \|\Pi_\lambda f \|_{L^p( \gamma)} \leq \varepsilon(\lambda)  \lambda^{\check{\delta}(p)}\|f\|_{L^2(M)},\qquad \text{ with } \lim_{\lambda\rightarrow +\infty}\varepsilon(\lambda) = 0. $$
 Writing a partition of unity
$$ 1= \sum_{n \in \mathbb{Z}} \Pi_{\frac n 2} $$
 and inserting another projector $\widetilde{\Pi}_\lambda=1_{\sqrt{- \mathbf{\Delta}} -\lambda \in [-1, 1]} $, 
$$\chi( \sqrt{- \Delta} -\lambda)= \sum_{n\in \mathbb{Z}}\Pi_{\frac n 2} \widetilde{\Pi}_{\frac n 2} \chi( \sqrt{- \Delta} -\lambda), $$
we obtain
$$ \|\chi( \sqrt{- \Delta} -\lambda) \|_{\mathcal{L}( L^2(M); L^p( \gamma))} \leq \sum_n \varepsilon(n) n^{\check{\delta}(p)} \|\widetilde{\Pi}_{\frac n 2} \chi( \sqrt{- \Delta} -\lambda) \|_{\mathcal{L}( L^2(M))}.$$
But due to the rapid decay of the function $\chi$, 
$$\|\widetilde{\Pi}_{\frac n 2} \chi( \sqrt{- \Delta} -\lambda) \|_{\mathcal{L}( L^2(M))}\leq \frac {C_N} {( 1+ |\lambda -\frac n 2|)^N}$$
and consequently we obtain
$$\|\chi( \sqrt{- \Delta} -\lambda) \|_{\mathcal{L}( L^2(M); L^p( \gamma))} \leq  o(1) \lambda^{\check{\delta}(p)}\qquad \text{ as }{\lambda \rightarrow + \infty}, $$ which is contradicting  our lower bounds for the spectral projector $\chi_{\lambda}$.\par
  
Now, if the manifold $M$ is a sphere, we know that $\Pi_\lambda$ has range either $\{0\}$ or the finite dimensional space or spherical harmonics of degree $k$ (corresponding to eigenvalues $-k(k+1)$) if $\sqrt{k(k+1)}-\lambda \in [0, \frac 1 2)$ (and in that case $k$ is uniquely determined). The sequence $(\lambda_n)$ appearing in~\eqref{eq.seq} corresponds necessarily to that latter case and the sequence $g_n=\Pi_{\lambda_n}f_n$ satisfies the  claimed lower bound.    
\end{proof}
\begin{remarque}Inequality~\eqref{eq.opt1} is optimized on the sphere by the so called {\em zonal} spherical harmonics which concentrate in $\mathcal{O}(\lambda^{-1})$ neighborhoods of two opposite points (the poles) and this can be easily seen by using asymptotic descriptions of these functions (see~\cite{Sz74}). Inequality~\eqref{eq.opt2} is optimized by the {\em highest weight} spherical harmonics $H_n(x_1, x_2, x_3) =  (x_1+ i x_2) ^n$ which concentrate on the equator $\{(x_1, x_2, x_3 ) \in \mathbb{S}^2; x_3=0\}$. Finally, inequality~\eqref{eq.opt3} is optimized {\em for $p=2$ and if the curve is a parallel circle} by some spherical harmonics which can be constructed by some {\em turning point} semi-classical analysis and concentrate on the parallel circles. In the general case (general curve or $2<p<4$), as suggested by one of the referees, one can check that the following spherical harmonics optimize inequality~\eqref{eq.opt3}. Let $H_{n, \phi}$ be the spherical harmonic rotated by angle $\phi$:
$$ H_{n, \phi}(x_1, x_2, x_3)= (x_1 + i (\cos (\phi) x_2 + \sin(\phi) x_3)) ^n$$
and for $\psi\in C^\infty_0(-1,1), \int \psi(s) ds =1$ and $\delta>0$ small enough
$$ u_n(x_1, x_2, x_3) = \int_{\phi} \Psi\Bigl( \frac{ N^{1/3} \phi}{ \delta}\Bigr) H_{n, \phi}(x_1, x_2, x_3)d \phi$$
(Remark that $u_n$ is clearly a spherical harmonic )
\end{remarque}
\section{The higher dimensional  case}
\subsection{The estimates in Theorem~\ref{th.variete}}
In this section we will give the modifications required to handle the case of submanifolds of higher dimensional manifolds. Following the approach in Section~\ref{sec.geod}, we define
for 
$$z\sim 0\in \mathbb{R}^{k}\rightarrow x(z)\sim 0 \in \mathbb{R}^{d}$$
a system of coordinates on $\Sigma$,
$$ T_\lambda (f)=\mathcal{T}_\lambda (f) (x(z))$$ and denote by $\mathcal{K}(x,x')$ the kernel of the operator $\mathcal{T}_\lambda \mathcal{T}^{\star}_\lambda$ and 
$$ K(z,z')= \mathcal{K} (x(z), x(z'))$$ the kernel of the operator $T_\lambda T^\star _\lambda$.  The starting point of our analysis is the higher dimensional analog of Lemma~\ref{geom}
\begin{lemme}\label{geombis}
There exist $\varepsilon\ll 1$, 
$(a_n^{\pm}, b_n )_{n \in \mathbb{N}}\in C^\infty(\R^d\times\R^d\times \R)$, 
such that for $\|x-x'\|\gtrsim \lambda^{-1}$ and any $N \in \mathbb{N}^*$,
\begin{equation}\label{eq.6.1}
{\mathcal K}(x,x')=\sum_{\pm}\sum_{n=0}^{N-1}\frac{e^{\pm i\lambda \psi(x,x')}}
{(\lambda\|x-x'\|)^{\frac{d-1}{2}+n}}\,
a_n^{\pm}(x,x', \lambda) + b_N (x, x', \lambda)\, ,
\end{equation}
where $\psi(x, x')$ is now the geodesic distance between the points $x$ and $x'$.
Moreover $a^{\pm}$ are real, have supports of size ${\mathcal O}(\varepsilon)$ with
respect to the first two variables  and are uniformly bounded with respect to $\lambda$. Finally 
$$ |b_N(x, x', \lambda)| \leq C ( \lambda|x-x'|)^{-\frac {d-1} 2 -N}.
$$
\end{lemme}

The proof of this Lemma is essentially the same as the proof of Lemma~\ref{geom} (see also~\cite[Lemma 2.7]{BuGeTz04}. We simply apply the stationnary phase formula in its full strength (i.e. with the asymptotic expansion to order $P$, see for example \cite[Theorem 7.7.5]{Ho83}), the two critical values of the phase ($\pm \psi(x,x')$) being evaluated as in Section~\ref{sec.courbe}.
\par
We shall also use the following Young's and Hardy-Littlewood inequalities
\begin{lemme}\label{lem.convolbis}
 Consider $2\leq p \leq+\infty$ and $T$ a convolution operator
$$ Tf(x)= \int_{\mathbb{R}^k} K(x-y) f(y) dy.$$
Then the norm of $T$ as an operator from $L^{p'}( \mathbb{R}^k)$ to $L^p(\mathbb{R}^k)$ is bounded by $\|K\|_{L^{\frac p 2}}$.
Furthermore, if $K(x)= \frac 1 {|x|^{2k/p}}$ and $p>2$, then the operator $T$ is still bounded from $L^{p'}( \mathbb{R}^k)$ to $L^p( \mathbb{R}^k)$ (despite the logarithmic divergence of $\|K\|_{L^{\frac p 2 }( \mathbb{R}^k)}$).
\end{lemme}
 
To conclude the proof of Theorem~\ref{th.variete} if $k<\frac {d-1}2$ is now easy: we simply remark that, since
$$d_g\left(x(z),x(z')\right)\sim \|z-z'\|,$$
  the kernel of $T_\lambda$, $K(x(z), x(z'))$ is bounded by 
\begin{equation}\label{eq.estim6}
 C (1+ \lambda \|z-z'\|)^{- \frac {d-1} 2}
\end{equation}
and apply the Young's inequality.\par
If $k< \frac{d-1} 2$, we conclude that the operator $T_\lambda T^\star _\lambda$ is bounded from $L^{2}(\Sigma)$ to $L^2( \Sigma)$ by $C \lambda^{- k}$ (resp. from $L^{1}(\Sigma)$ to $L^\infty( \Sigma)$ by $C$) and consequently the operator $T_\lambda$ is bounded from $L^2 ( M)$ to $L^2( \Sigma)$ by  $C \lambda^{- k/2}$ (resp. from $L^2 ( M)$ to $L^\infty( \Sigma)$ by  $C $). Taking into account that
\begin{equation}
\label{eq.6.2}
 \varphi_\lambda = \lambda^{\frac{d-1} 2} \mathcal{T}_\lambda (\varphi_\lambda)
\end{equation}
we obtain (if $k< \frac{d-1} 2 $) Theorem~\ref{th.variete} for $p=2$ and $p= + \infty$  and consequently, by interpolation, for any $2\leq p \leq + \infty$. \par
If $k\geq \frac{d-1} 2$, the simple argument above still applies (using Hardy Littlewood inequality), but it gives Theorem~\ref{th.variete} only in the range 
\begin{equation}
\begin{cases}  p\geq \frac{4k} {d-1} &\text{ if } k>\frac{ d-1} 2,\\
p>\frac{4k} {d-1} &\text{ if } k= \frac {d-1} 2.
\end{cases} 
\end{equation}
 To obtain the full result we have to refine the analysis and take advantage of some oscillations in the phase in~\eqref{eq.6.1}. 
We first consider the contribution of $\{(z,z'); |z-z'|\leq C \lambda^{-1}\}$ and define for $\chi\in C^\infty_0( \mathbb{R}^d)$
$$ K_0(z,z)=K(z,z)\chi(\lambda(z-z')).$$
Then, using that $K_0$ is bounded and supported in the set $\{(z,z'); \|z-z'\|\leq C \lambda^{-1}\}$,  in view of Young's inequality, the contribution of $K_0$ is easily dealt with. 
 Furthermore, the analysis above is sufficient, if $N$ is chosen large enough in Lemma~\ref{geombis}, to deal with the contribution of the remainder term $b_N$.
Indeed, if $N> \frac{d-1} 2$
$$ \int_{ \lambda^{-1}\leq |z| \leq 1} \frac{dz} {(\lambda |z|)^{\frac {d-1} 2 +N}}\leq C\lambda^{-k}.$$
As a consequence, we will focus on the contributions of the main term $a_0^+$ and abusing notations we will still denote by $\mathcal{K}$ the function 
$$\mathcal{K}(x,x')=\frac{e^{ i\lambda \psi(x,x')}}
{(\lambda\|x-x'\|)^{\frac{d-1}{2}}} a_0^{+}(x,x', \lambda)
$$
 The contributions of the other terms ($a_0^-,a_n^\pm, 1\leq n <N$) could be dealt with similarly, as will appear clearly in the proof below.
  
We now consider a partition of unity on $B= \{x\in \mathbb{R}^{k}; |x|<1\}$,
$$ 1= \chi(\lambda x) +\sum_{j=1}^{\log\lambda/\log 2} \widetilde\chi (2^j x)$$
where $\widetilde \chi\in C^\infty_0( \mathbb{R}^d)$ and $\widetilde \chi$ is supported in the set $\{x; \frac 1 2 < |x| < 2\}$. We consider the related partition of the kernel:
\begin{multline}
 K(z,z)=K(z,z)\chi(\lambda(z-z')) +\sum_{j=1}^{\log\lambda/\log 2}K(z,z') \widetilde\chi (2^j (z-z'))\\
\equiv K_0(z, z')+ \sum_{j=1}^{\log\lambda/\log 2}K_j(z,z')
\end{multline}

The main step in the proof of Theorem~\ref{th.variete} is the proof of
\begin{proposition}\label{prop.base}
If $j$ is large enough, the operator $(TT^\star)_j$ whose kernel is $K_j(z,z')$ satisfies the bounds
\begin{equation} 
\begin{aligned}
\|(TT^\star)_jf\|_{L^\infty(\Sigma)} &\leq C \left(\frac{ 2^j} \lambda \right)^{\frac {d-1} 2}\|f\|_{L^1( \Sigma)}\\
\|(TT^\star)_jf\|_{L^2(\Sigma)} &\leq C 2^{-jk}\left(\frac { 2^j} \lambda\right)^{\frac {d-1} 2 + \frac{k-1} 2}\|f\|_{L^2( \Sigma)} 
\end{aligned} 
\end{equation}
\end{proposition}
Let us first show how Proposition~\ref{prop.base} implies Theorem~\ref{th.variete}:       by interpolation, we deduce
$$\|(TT^\star)_jf\|_{L^p(\Sigma)} \leq C 2^{-\frac{2jk} p}\left(\frac { 2^j} \lambda\right)^{\frac {d-1} 2 + \frac{k-1} p} \|f\|_{L^{p'}( \Sigma)}$$
As a consequence, summing for $j \leq \log(\lambda) /\log(2)$ we obtain
\begin{multline}
\|TT^\star f\|_{\mathcal{L}(L^{p'}(\Sigma); L^p( \Sigma))} \leq C \lambda^{-\frac{d-1} 2 -\frac {k-1} p}\sum_{j \leq \log( \lambda) /\log(2)}  2^{j( \frac {d-1} 2 - \frac {k+1} p)}\\
\leq \begin{cases}C \lambda^{-\frac{d-1} 2 - \frac {k-1} p+ \frac {d-1} 2 - \frac {k+1} p} & \text{ if } p > \frac {2(k+1) } {d-1}\\
C \lambda^{-\frac{d-1} 2 - \frac {k-1} p}& \text{ if } p < \frac {2(k+1) } {d-1}\\
C \lambda^{-\frac{d-1} 2 - \frac {k-1} p}\log( \lambda) & \text{ if } p = \frac {2(k+1) } {d-1}
\end{cases}
\end{multline}
 which, taking~\eqref{eq.6.2} into account gives the estimates in Theorem~\ref{th.variete} (remark that the case $p<2( k+1) /(d-1)$ occurs only when $k=d-1$). Remark also that the fact that we can only prove the estimate for large $j$ is not a problem, as, shrinking the support of the symbol $a$ (i.e. taking $\epsilon$ small enough)  in Theorem~\ref{sogge}, the contributions of small $j$ vanish.
\par
We now come back to the proof of Proposition~\ref{prop.base}.
The $L^1-L^\infty$ bound is straightforward and we can focus on the $L^2$ bound.
We have
$$ K_j(z,z')= e^{i \lambda d(z,z')} \frac {\widetilde \chi( 2^j( z-z'))}{(\lambda d(z,z'))^{\frac {d-1} 2}} a_0 (z,z')$$
where $d(z,z')$ is the geodesic distance (in $M$) between points $x(z)$ and $x(z')$. As a consequence, on the support of the r.h.s., $d(z,z') \sim 2^{-j}$.
Next we introduce a partition of unity locally finite (uniformly with respect to $j$):
$$ 1= \sum_{p \in \mathbb{Z}^{d-1}} \chi( 2^j z-p),  $$
and write
$$ K_j(z,z') = \sum_{q,\widetilde q \in \mathbb{Z}^{d-1}}\chi( 2^jz-q) K_j(z,z') \chi( 2^jz'-\widetilde q)$$
We denote by $R_{j,q,\widetilde q}$ the operator whose kernel is 
$$\chi( 2^jz-q) K_j(z,z') \chi( 2^jz'-\widetilde q).$$
Remark that due to the support properties of $K_j$, in the expression above, only the contributions from $(q,\widetilde q)$ such that $\|q-\widetilde q\|\leq C$ do not vanish. Remark also that by quasi-orthogonality in $L^2$ (due to the fact that the partition of unity is locally finite), we have 
$$\|(TT^\star)_j\|_{ \mathcal{L} ( L^2( \Sigma))} \leq C \sup_{q,\widetilde q} \| R_{j,q,\widetilde q}\|_{ \mathcal{L} ( L^2( \Sigma))}$$
For simplicity, we shall only estimate the norm of $  R_{j,q,q}$ (the case of $(q,\widetilde q), \|q-\widetilde q\|\leq C$ is similar). Using a translation and an orthogonal (linear) transformation, we can assume that $q=0$ and the metric satisfies $g_{m,n}(q)= Id$. We now perform a change of variables, set $Z= 2^j z$ and obtain
$$ \| R_{j,0,0}\|_{ \mathcal{L} ( L^2( \Sigma))}= 2^{-j k}\| \underline{ R_{j,0,0}}\|_{ \mathcal{L} ( L^2( \Sigma))}$$
where the kernel of the operator $ \underline {R_{j,0,0}}$, $\underline{K}(Z, Z')$ satisfies
\begin{multline}
 \underline{K} (Z, Z')= \chi(Z) \widetilde{\chi}( Z')K_j(2^{-j} z, 2^{-j} z')\\
 = e^{i \lambda d(2^{-j}Z,2^{-j} Z')} \chi(Z) \chi(Z')
\frac{\widetilde\chi(Z-Z')} 
{ (\lambda d(2^{-j}Z,2^{-j} Z'))^{\frac {d-1} 2}} a_0 (2^{-j}Z,2^{-j} Z')\\
=\left( \frac {2^j} { \lambda}\right)^{\frac{d-1} 2} e^{i \lambda2^{-j} d_j(Z, Z')} \chi(Z) \chi(Z')\frac{\widetilde\chi(Z-Z')} { d_j(Z, Z')^{\frac {d-1} 2}} a_0 (2^{-j}Z,2^{-j} Z')
\end{multline}
where $d_j(Z,Z')= 2^j d(2^{-j}Z, 2^{-j} Z') $ is the distance (measured in $\mathbb{R}^d$) between the points $(Z,0)$ and $(Z',0)$ for the family of metrics $g_{m,n}^j (X)= g_{m,n}(2^{-j}X)$ which converges as $j$ tend to the infinity to the metric $Id$ (in the $C^\infty$ topology). As a consequence (and using that on the support of the r.h.s., $\frac 1 2\leq \|Z-Z'\| \leq 2$), our kernel has the following form
$$\underline{K} (Z, Z')=\left( \frac {2^j} { \lambda}\right)^{\frac{d-1} 2} e^{i \lambda2^{-j} d_j(Z, Z')} \sigma(Z,Z',j)$$
with $d_j(Z,Z')$ arbitrarily close (for large $j$) in $C^\infty$ topology to $\|Z-Z'\|$ and $\sigma$ a function uniformly bounded with respect to $j$ in $C^\infty_0$ topology. Proposition~\ref{prop.base} is now a consequence of the following non degeneracy property:
\begin{proposition}\label{prop.ter}
Consider an operator $T$ on $L^2( \mathbb{R}^k)$, whose kernel $K(Z,Z')$ has the following form
$${K} (Z, Z')= e^{i \mu d_j(Z, Z')} \sigma(Z,Z',j)$$
with $d_j(Z,Z')$ arbitrarily close (for large $j$) in $C^\infty$ topology to $d_\infty(Z,Z')=\|Z-Z'\|$ and $\sigma$ a function uniformly bounded with respect to $j$ in $C^\infty$ topology supported in the set 
$$\{(Z,Z'); \|Z\|\leq 1, \frac 1 2 \leq \|Z-Z'\|\leq 2\}.$$
 Then, for large $j$, it satisfies the bound 
$$
\|T\|_{\mathcal{L} (L^2( \mathbb{R}^k))} \leq C \mu^{-\frac {k-1} 2}.$$
\end{proposition}
The proof of this result (for $j = \infty$) is standard. We shall recall it to check that it goes through with the parameter $j$. By using partition of unities, we can assume that $\sigma$ is supported in the set 
$$\{(Z,Z'); \|Z-Z_0\|\leq \epsilon, \|Z'-Z'_0\|\leq \epsilon, \|Z_0-Z'_0\|\in [\frac 1 2, 2]\}$$
We take polar coordinates centered in $Z_0$ and write $Z'=Z_0+ r\theta$. We can assume that $Z'_0= (r_0, 0,0)$. 
We have 
$$\nabla_Z(\|Z- Z'\|)\mid_{Z=Z_0}= - \theta.$$ 
As a consequence, taking as coordinates $Z_1, X= (Z_2,Z_3)$, we have 
$$ \nabla^2_{X, \theta}( \|Z-Z'\|)\mid_{Z=Z_0, Z'= Z'_0}= - \text{Id}$$
and we can apply to the operator $T_{Z_1,r}$ obtained by freezing the $T$ and $r$ variables the following classical non-degenerate phase property (see for example~\cite[Proposition IX.1.1]{St93})

\begin{lemme}\label{lem.nondegene} Consider an operator 
$$ T_\lambda f( \xi) = \int_{\mathbb{R}^n} e^{i \lambda \Phi(x,  \xi)} \Psi (x, \xi) f(x) dx$$
where $\Psi(x, \xi)$ is a fixed smooth function of compact support in $\mathbb{R}^n$ and the phase $\Phi$ is real valued, smooth and satisfies
$$ \text{det} \left( \frac{ \partial^2 \Phi(x, \xi)} {\partial x_i \partial\xi_j}\right) \neq 0.$$ Then
$$\|T_\lambda f \|_{L^2( \mathbb{R}^n)} \leq C\lambda^{-n/2} \|f\|_{L^2( \mathbb{R}^n)}.$$
\end{lemme}
A simple use of Minkovski inequality gives Proposition~\ref{prop.ter} (in the case $j= + \infty$). To conclude the proof for large $j$, we remark that Lemma~\ref{lem.nondegene} is stable by small (smooth) perturbations.
\subsection{Optimality of Theorem~\ref{th.variete}}
The first regime to take into account for the optimality is the zonal regime. If we consider functions on $\mathbb{S}^d$ depending only on the geodesic
distance to a fixed point, we obtain the zonal eigenfunctions on  $\mathbb{S}^d$. The zonal eigenfunctions can be expressed in terms of zonal
spherical harmonics which in their turn can be expressed in terms of the classical Jacobi polynomials (see e.g. \cite{So93}).
 In that case we can show that we have a pointwise concentration. If $Z_n$ is the $n$-th zonal eigenfunction (with eigenvalue $\lambda^2 = n(n+d-1)$)  
\begin{equation}\label{point}
|Z_{n}(x)|\approx  n^{\frac{d-1}{2}}\|Z_n\|_{L^2(\mathbb{S}^d)},\quad d(x, x_0)\leq \frac cn
\end{equation}
As a consequence, for any submanifold $\Sigma\subset \mathbb{S}^d$ if we choose a pole $P \in \Sigma$, and consider the family of corresponding zonal eigenfunctions, $Z_n$, we obtain
$$ \|Z_n\mid_\Sigma\|_{L^p( \Sigma)}\geq c n^{\frac{d-1} 2 - \frac k p}\|Z_n\|_{L^2(\mathbb{S}^d)}, \qquad c>0.$$
and this shows the optimality of Theorem~\ref{th.variete} if $k \leq d-2$ and if $k=d-1$ and $\frac {2d} {d-1} \leq p \leq + \infty$.\par
To obtain the optimality of Theorem~\ref{th.variete} in the last regime ($k=d-1$, $2 \leq p \leq \frac{2d} {d-1}$), we turn to the highest weight eigenfunctions $e_n=n^{\frac {d-1} 4}(x_1+ i x_2)^n$ (corresponding to eigenvalues $\lambda^2= n(n+d-1)$). In that case, if $\Sigma$ contains the geodesic $$\gamma = \{ x= (x_1, \dots, x_d); x_3= \dots =x_d =0\},$$ then we obtain
$$ \|e_n \mid_{\Sigma}\|_{L^p( \Sigma)} \sim n^{\frac{d-1} 4 - \frac{k-1} {2p}}\|e_n\|_{L^2(\mathbb{S}^d)}.$$

\def\cprime{$'$}

\end{document}